\documentclass{amsart}

\usepackage{amssymb,amscd,amsthm,latexsym, amsmath}
\usepackage[all]{xy}
\usepackage{color}

\newtheorem{defi}{Definition}[section]
\newtheorem{theo}{Theorem}[section]
\newtheorem{lemma}{Lemma}[section]
\newtheorem{prop}{Proposition}[section]
\newtheorem{coro}{Corollary}[section]

\def\into{ \rightarrowtail }

\def\splito{ \rightleftarrows }

\newdir{ >}{{}*!/-8pt/@{>}}

\def\EE{ \mathbb{E} }

\def\RR{ \mathbb{R} }
\def\GG{ \mathbb{G} }
\def\CC{ \mathbb{C} }
\def\DD{ \mathbb{D} }

\def\Pt{ \mathrm{Pt} }

\newcommand{\Gp}{\mathsf{Gp}}

\newcommand{\Hom}{\mathsf{Hom}}

\newcommand{\TC}{\mathsf{T\sharp C}}
\newcommand{\TXC}{\mathsf{T_{X_{\bullet}}\sharp C}}
\newcommand{\TXTC}{\mathsf{T_{X^T_{\bullet}}\sharp C}}
\newcommand{\GGC}{\mathsf{\GG\sharp C}}
\newcommand{\ColgC}{\mathsf{ColgC}}

\newcommand{\Alg}{\mathsf{Alg}}
\newcommand{\AlgT}{\mathsf{AlgT}}

\begin{document}

\title{Kleisli categories, $T$-categories and internal categories}

\author{Dominique Bourn}
\address{Univ. Littoral C\^ote d'Opale, UR 2597, LMPA,
	Laboratoire de Math\'ematiques Pures et Appliqu\'ees Joseph Liouville,
	F-62100 Calais, France}
\email{bourn@univ-littoral.fr}
\keywords{Internal categories and groupoids, Kleisli categories, $T$-categories, internal $n$-categories and $n$-groupoids, operads, multicategories. \\{\small 2020 {\it Mathematics Subject Classification.} Primary 18C20,18D40,18M60,18M65,18N50.}}

\date{}

\begin{abstract}
	We investigate the properties of the Kleisli category $KlT$ of a monad $(T,\lambda,\mu)$ on a category $\EE$ and in particular the existence of (some kind of) pullbacks. This culminates when the monad is cartesian. In this case, we show that any $T$-category in $\EE$ in the sense of A. Burroni coincides with a special kind of internal category in $KlT$. So, it is the case in particular for $T$-operads and $T$-multicategories. More unexpectedly, this, in turn, sheds new lights on internal categories and $n$-categories.
\end{abstract}

\maketitle

\section*{Introduction}

The notion of $T$-category in a category $\EE$ endowed with a monad $(T,\lambda,\mu)$ was introduced by A. Burroni in \cite{AB}. It is a kind of a mix of a relational algebra on $T$ in the sense of Barr \cite{Barr1} and of "something" which looks like an internal category, but shifted by this monad, see precise definition in Section \ref{Tcat}. The aim of this work is to investigate what is exactly this "something".
It easily appears that a $T$-category is a special kind of $3$-truncated simplicial object in the Kleisli category $KlT$ of this monad. Now, when, in addition, the monad $(T,\lambda,\mu)$ is cartesian, we shall show that:\\
1) the bijective on object natural functor $\bar F^T: \EE \to KlT$ is actually an inclusion;\\
2) the subcategory $\EE$ then appears to be \emph{left cancelable} in $KlT$, i.e. such that $h\in \EE$ and $g\in \EE$ imply $f\in\EE$ when $h=g.f$ in $KlT$;\\
3) when $\EE$ is finitely complete, not only the category $\AlgT$ of algebras on the monad is finitely complete as well, but the Kleisli category $KlT$, which is not finitely complete, is however such that any map $f: X\to Y$ in (the subcategory) $\EE$ has a pullback along any map in $KlT$ which still belongs to $\EE$.

From this observation, we shall show that, when the monad is cartesian, the previous $3$-truncated simplical object in $KlT$ is actually underlying a regular internal category in $KlT$, and that any $T$-category coincides with this kind of internal category. So, according to \cite{Lein}, $T$-operads and $T$-multicategories appear to be internal categories in $KlT$.

If the properties of the category $\AlgT$ of $T$-algebras are well known, those of $KlT$ have been neglected. This work gave us the opportunity to investigate how, step by step, the  several assumptions of a cartesian monad surprisingly organize the properties of $KlT$.

Conversely, and more unexpectedly, any internal category in $\EE$ will appear to be a special kind of $\GG$-category where $(\GG,\sigma,\pi)$ is the monad on the category $\Pt\EE$ of split epimorphisms in $\EE$ whose category of algebras was known to be nothing but the category $Grd\EE$ of internal groupoids  in $\EE$, see \cite{B-1}.
From that we shall how internal $n$-categories and $n$-groupoids are related to this monad as well. So that  the following whole tower of fibrations is entirely ruled by the monad $(\GG,\sigma,\pi)$ and constructed from it:
$$...\; n{\rm -}Cat\EE \stackrel{(\;)_{n-1}}{\rightarrow} (n-1){\rm -}Cat\EE \; ......\; 2{\rm -}Cat\EE \stackrel{(\;)_1}{\rightarrow} Cat\EE \stackrel{(\;)_0}{\rightarrow} \EE$$
Beyong the heuristic interest of this result, it could appear very useful when combinatorial diagrammatic calculations will have been developed on computers. For some fresh results on Burroni's $T$-categories in another direction see \cite{Th}.

The ideas of this work came to my mind during a talk of M. Batanin in Nice for the \emph{Homotopical days} (7-9 Dec. 2022), see \cite{Bat2} and also \cite{Bat1}. The article is organized along the following lines:\\ 
Section 1: brief recalls about monads. Section 2: properties of the Kleisli category $KlT$. Section 3: brief recalls about internal categories. Section 4: recalls about internal groupoids and the monad $(\GG,\sigma,\pi)$.  Section 5: recalls about $T$-categories. Section 6: when $T$-categories in $\EE$ coincide with a special kind of internal categories in $KlT$. Section 7: when internal categories in $\EE$ coincide with a special kind of $\GG$-categories. Section 8: extensions of the results of the previous section to internal $n$-categories and $n$-groupoids. Any internal category $X_{\bullet}$ produces a cartesian monad $(T_{X_{\bullet}},\lambda_{X_{\bullet}},\mu_{X_{\bullet}})$ on the slice category $\EE/X_0$; Section 9 is devoted to make explicit all the results of Section 6 about this monad. Section 10 is devoted to translate  the results of this same section to $T$-operads, $T$-multicategories, and their algebras.

\section{Monads}

\subsection{Basics}

Let us briefly recall the basics on monads. A monad on a category $\EE$ is a triple $(T,\lambda,\mu)$ of an endofunctor $T$ and two natural transformations: $$Id_{\EE} \stackrel{\lambda}{\longrightarrow} T \stackrel{\mu}{\longleftarrow} T^{2}$$
satisfying $\mu.\mu T=\mu.T\mu$ and $\mu.\lambda T=1_{T}=\mu.T\lambda$. An adjoint pair $(U,F,\lambda,\epsilon): \bar\EE \splito \EE$ determines  the monad $(T,\lambda,\mu)=(U.F,\lambda,U\epsilon F)$ on $\EE$.

A $T$-\emph{algebra} \cite{EM} on an object $X$ is given by a map $\xi: T(X)\to X$ satisfying $\xi.\lambda_X=1_X$ and $\xi.\mu_X=\xi.T(\xi)$. Accordingly the pair $(T(X),\mu_X)$ produces a $T$-algebra on $T(X)$. A morphism $f: (X,\xi) \to (Y,\gamma)$ of $T$-algebras is given by a map $f: X \to Y$ such that $f.\xi=\gamma.T(f)$. 

This construction determines the category $\AlgT$ of $T$-algebras and the forgetful functor $U^T: \AlgT \to \EE: (X,\xi) \mapsto X$ which is obviously conservative. It has the functor $F^T: \EE \to \AlgT$ defined by $F^T(X)=(T(X),\mu_X)$ as left adjoint which makes $U^T$ a left exact functor. The monad associated with the adjoint pair  $(U^T,F^T)$ recovers the initial monad $(T,\lambda,\mu)$. From an adjoint pair $(U,F,\lambda,\epsilon): \bar\EE \splito \EE$ and its associated monad $(T,\lambda,\mu)=(U.F,\lambda,U\epsilon F)$ we get a comparison functor $A_{(U,F)}: \bar\EE \to \AlgT$, defined by $A_{(U,F)}(Z)=(U(Z),U(\epsilon_Z))$,  making the following adjoint pairs commute:
$$\xymatrix@=18pt{
	\bar\EE  \ar[rr]^{A_{(U,F)}} \ar[dr]_F  &&	\AlgT \ar[dl]_{F^T} \\
	& \EE  \ar@<1ex>[ur]_{U^T}  \ar@<1ex>[ul]_U  &
	 }
$$
The functor $U$ is said to be \emph{monadic} when the  comparison functor $A_{(U,F)}$ is an equivalence of categories.
 
A comonad $(C,\epsilon,\nu)$ is the dual of a monad; it determines the category $\ColgC$ of co-algebras and a coadjoint pair $(F_C,U_C): \ColgC \splito \EE$.
Any adjoint pair $(U,F,\lambda,\epsilon): \bar\EE \splito \EE$ determines the comonad $(C,\epsilon,\nu)=(F.U,\epsilon,F\lambda U)$ on $\bar\EE$ and a comparison functor $C_{(U,F)}: \EE\to \ColgC$,  defined by $C_{(U,F)}(X)=(F(X),F(\lambda_{X}))$, making the following coadjoint pairs commute:
$$\xymatrix@=18pt{
	\EE  \ar[rr]^{C_{(U,F)}} \ar[dr]_U  &&	\ColgC \ar[dl]_{U_C} \\
	& \bar\EE  \ar@<1ex>[ur]_{F_C}  \ar@<1ex>[ul]_F  &
}
$$
The functor $F$ is said to be \emph{comonadic} when the  comparison functor $C_{(U,F)}$ is an equivalence of categories.

\subsection{Cartesian monads}

A functor $F:\CC\to\DD$ is said to be \emph{cartesian} when $\CC$ has pullbacks and $F$ preserves them. A natural transformation $\nu: F\Rightarrow G$ between any pair of functors is said to be \emph{cartesian} when, given any map $f: X\to Y\in \CC$, the following square is a pullback in $\DD$:
$$
\xymatrix@=4pt{
	F(X) \ar[dd]_{F(f)}  \ar[rr]^{\nu_X}  && G(X) \ar[dd]^{G(f)}\\
	&&&&   \\
	F(Y) \ar[rr]_{\nu_Y}  &&  G(Y)
}
$$
A monad $(T,\lambda,\mu)$ is \emph{cartesian} when the three ingredients are cartesian. 
A. Burroni \cite{AB} was deeply involved in the cartesian monad induced by the free adjunction $Cat\splito Gph$ between categories and directed graphs. More specifically, from T. Leinster \cite{Lein2}, the free monoid monad $(M,\lambda,\mu)$ on $Set$ is a cartesian one: it is the restriction of Burroni's monad to directed graphs with only one object.
Sections \ref{T_X} will be devoted to a cartesian monad associated with any internal category $X_{\bullet}$. 
\begin{prop}\label{lcart}
	Let $(T,\lambda,\mu)$ be a monad on $\EE$ where $\mu$ is cartesian; the two conditions are equivalent:\\
	1) $\lambda$ is the equalizer of $\lambda T$ and $T\lambda$, and
	2) $\lambda$ is cartesian. 
\end{prop}
\proof
When $\mu$ is cartesian, the natural transformations $\lambda T$ and $T\lambda$ are necessarily cartesian as well, being splittings of the cartesian $\mu$. 
Now, for any map $f: X\to Y$ in $\Sigma$, consider the following diagram:
$$\xymatrix@=7pt{
	X \ar[dd]_{f}  \ar@{ >->}[rr]^{\lambda_X}  && T(X) \ar[dd]_{T(f)}  \ar@<1ex>[rr]^{\lambda_{T(X)}} \ar@<-1ex>[rr]_{T(\lambda_{X})} && T^2(X) \ar[dd]^{T^2(f)} \ar[ll]|{\mu_X}\\
	&&&&   \\
	Y \ar@{ >->}[rr]_{\lambda_Y}  &&  T(Y) \ar@<1ex>[rr]^{\lambda_{T(Y)}} \ar@<-1ex>[rr]_{T(\lambda_{Y})} && T^2(Y) \ar[ll]|{\mu_Y}
}
$$
Any of the right hand side commutative squares is a pullback.
Moreover, under  assumption 1), the two horizontal ones  are pullbacks. Accordingly, the "box lemma" for pullbacks makes the left hand square a pullback as well.

In any case, by the following lemma, when $\lambda$ is cartesian, $\lambda$ is the equalizer of $\lambda T$ and $T\lambda$.
\endproof
\begin{lemma}\label{equali}
	Given any cosplit parallel pair in a category $\EE$:
	$$\xymatrix@=10pt{
		X \ar@<2ex>[rr]^{m} \ar@<-2ex>[rr]_{m'} && \check X \ar[ll]|g
	}
	$$
	any equality $m.k=m'.h$ implies $k=h$. Accordingly the pullback of the maps $m$ and $m'$ produces their equalizer. So:\\
	1) given any monad $(T,\lambda,\mu)$, if $\lambda$ is cartesian, then $\lambda$ is the equalizer of the pair $(\lambda_{T},T(\lambda))$;\\
	2) any cartesian functor preserves the equalizers of cosplit parallel pairs. 
\end{lemma}
Later on, we shall need the following:
\begin{defi}
	A monad $(T,\lambda,\mu)$ is said to be half-cartesian when the endofunctor $T$ is cartesian and $\lambda$ is the equalizer of the pair $(\lambda_{T},T(\lambda))$.
\end{defi} 
Accordingly, a monad is cartesian if and only if it is half-cartesian and $\mu$ is cartesian.

\subsection{Cartesian adjoint pairs}

It is then natural to call \emph{cartesian adjoint pair}, any adjoint pair $(U,F):\CC\splito \DD$ such that $\CC$ and $\DD$ has pullbacks, the functor $F$ is cartesian, the natural transformations $\lambda: Id_\DD\Rightarrow U.F$ and $\epsilon: F.U \Rightarrow Id_{\CC}$ are cartesian. Then the induced monad $(T,\lambda,\mu)$ on $\DD$ is clearly a cartesian monad, since $\mu=U(\epsilon_F)$. The induced comonad $(C,\epsilon,\nu)$ on $\CC$ is cartesian as well: the functor $C=F.U$ is cartesian and the natural transformation $\nu$, being a section of the natural transformation $\epsilon_C$, is cartesian as soon as so is $\epsilon$.

\begin{prop}\label{cartadj}
	Given any cartesian adjoint pair $(U,F):\CC\splito \DD$, the natural transformation $\mu: T^2\to T$ of the induced monad on $\DD$ is such that the following  diagram is a kernel equivalence relation:
	$$\xymatrix@=5pt{
	T^3(X) \ar@<-2ex>[rrrr]_{\mu_{T(X)}} \ar@<2ex>[rrrr]^{T(\mu_{X})} &&&& T^2(X)\ar[rrrr]_{\mu_X}\ar[llll]|{T(\lambda_{T(X)})}&&&& T(X) 
	}
	$$
	Conversely suppose that $(T,\lambda,\mu)$ is a cartesian monad. The adjoint pair $(U^T,F^T): \Alg T\splito \EE$ is a cartesian one if and only if  the natural transformation $\mu$ satisfies the above property. In this case, given any $T$-algebra $\xi: T(X)\to X$, the following diagram produces a kernel equivalence relation:
	$$\xymatrix@=5pt{
		T^2(X) \ar@<-2ex>[rrrr]_{\mu_{X}} \ar@<2ex>[rrrr]^{T(\xi)} &&&& T(X)\ar[rrrr]_{\xi}\ar[llll]|{T(\lambda_{X})}&&&& X 
	}
	$$
\end{prop}
\proof
The first assertion is the consequence of the fact that the commutative square underlying the diagram in question is the image by the cartesian functor $U$ of the following pullback:
$$\xymatrix@=3pt{
	(F.U)^2(F(X))	\ar[dd]_{F.U(\epsilon_{F(X)})} \ar[rrrr]^{\epsilon_{F.U.F(X)}} &&&& F.U.F(X)   \ar[dd]^{\epsilon_{F(X)}}\\
	&&&&\\
	F.U.F(X) \ar[rrrr]_{\epsilon_{F(X)}} &&&& F(X) 
}
$$
Now suppose the monad is cartesian. Given any $T$-algebra $x: T(X)\to X$ on the object $X$ ans applying $\mu$-cartesianness to the map $x$, we get a pullback, in such a way that the map $T(\mu_{X})$ delineates the composition map of an internal category in $\EE$ (see Section \ref{ccategory} below):
$$\xymatrix@=20pt{
	T^3(X)\ar@<-2ex>[rr]_{T^2(x)} \ar[rr]|{T(\mu_{X})} \ar@<2ex>[rr]^{\mu_{T(X)}}  &&	T^2(X) \ar@<-2ex>[rr]_{T(x)} 
	\ar@<2ex>[rr]^{\mu_{X}}  && T(X)\ar[ll]|{T(\lambda_{X})}
}
$$
When, in addition, $\mu$ satisfies the condition in question, the pair $(\mu_{T(X)},T(\mu_{X}))$ is the kernel pair of $\mu_{X}$
and this category is actually a groupoid  (see Theorem \ref{carac} below).
Any morphism of $T$-algebras: $f: (X,x)\to (Y,y)$ produces the following vertical discrete fibration between groupoids:
$$\xymatrix@=15pt{
	T^3(X) \ar[dd]_{T^3(f)}\ar@<-2ex>[rr]_{T^2(x)} \ar[rr]|{T(\mu_{X})} \ar@<2ex>[rr]^{\mu_{T(X)}}  &&	T^2(X) \ar[dd]^{T^2(f)}\ar@<-2ex>[rr]_{T(x)} 
	\ar@<2ex>[rr]^{\mu_{X}}  && T(X)\ar[ll]|{T(\lambda_{X})} \ar[dd]^{T(f)}\\
	&&&&\\
	T^3(Y)\ar@<-2ex>[rr]_{T^2(y)} \ar[rr]|{T(\mu_Y)} \ar@<2ex>[rr]^{\mu_{T(Y)}}  &&	T^2(Y) \ar@<-2ex>[rr]_{T(y)} 
	\ar@<2ex>[rr]^{\mu_Y}  && T(Y)\ar[ll]|{T(\lambda_Y)}
}
$$
So, it is a discrete cofibration as well, and the following rightward left hand side commutative square is a pullback:
$$\xymatrix@=12pt{
	T^2(X) \ar[dd]_{T^2(f)} \ar@<-2ex>[rr]_{T(x)}   &&	T(X) \ar[dd]^{T(f)} 
	\ar@<-2ex>[rr]_{x} \ar[ll]_{\lambda_{T(X)}} && X\ar[ll]_{\lambda_{X}} \ar[dd]^{f}\\
	&&&&\\
	T^2(Y)\ar@<-2ex>[rr]_{T(y)}   &&	T(Y) \ar@<-2ex>[rr]_{y} \ar[ll]_{\lambda_{T(Y)}}
	  && Y\ar[ll]_{\lambda_Y}
}
$$
So, the rightward right hand side commutative square is a pullback as well, since, composed with the leftward right hand side pullback ($\lambda$ is cartesian), it gives rise to the pullback obtained by composition of the two left hand side pullbacks. This exactly means that the co-unit $\epsilon: F^T.U^T\rightarrow Id_{\AlgT}$ of the comonad on $\AlgT$ is cartesian. The last assertion is obtained by  applying $\epsilon$-cartesianness to $\epsilon$ itself.
\endproof
So, let us introduce the following:
\begin{defi}
	A monad $(T,\lambda,\mu)$ is said to be hypercartesian when it is cartesian and the natural transformation $\mu$ is such that:
	$$\xymatrix@=5pt{
	T^3(X) \ar@<-2ex>[rrrr]_{\mu_{T(X)}} \ar@<2ex>[rrrr]^{T(\mu_{X})} &&&& T^2(X)\ar[rrrr]_{\mu_X}\ar[llll]|{T(\lambda_{T(X)})}&&&& T(X) 
	}
	$$
	is a kernel equivalence relation for any object $X$.
\end{defi}

\subsection{Autonomous adjoint pairs}

Let us introduce the following:
\begin{defi}
	An adjoint pair $(U,F):\bar\EE\splito \EE$ is said to be autonomous when $U$ is monadic and $F$ comonadic.
\end{defi}
In other words, an adjunction is autonomous when it does not expand in new adjunctions via the algebra or co-algebra constructions.
The aim of this section is to prove that any half-cartesian monad makes the adjoint pair $(U^T,F^T): \AlgT\splito \EE$ an autonomous one.

\begin{prop}\label{Tff}
	Let $(U,F): \bar\EE \to \EE$ be  an adjoint pair and $(T,\lambda,\mu)$ its associated monad on $\EE$. The two following conditions are equivalent:\\
	1) the natural transformation $\lambda$ is the equalizer of the pair $(\lambda_{T},T(\lambda))$;\\
	2) the comparison functor $C_{(U,F)}:\EE\to \ColgC$ is fully faithful.\\
	Any of these conditions implies that the functor $F: \EE\to \bar\EE$ is conservative.
\end{prop}
\proof
Suppose 1). Let $h: F(X)\to F(Y)$ be a map in $\bar\EE$ making the following left hand side square commute:
$$\xymatrix@=7pt{
	F(X) \ar[dd]_{h}  \ar[rr]^{F(\lambda_X)}  && FT(X) \ar[dd]^{FU(h)}	&& X \ar@{.>}[dd]_{k}  \ar@{ >->}[rr]^{\lambda_X}  && T(X) \ar[dd]_{U(h)} \ar@<1ex>[rr]^{\lambda_{T(X)}} \ar@<-1ex>[rr]_{T(\lambda_{X})} && T^2(X) \ar[dd]^{TU(h)}\\
	&&&&&&   \\
	F(Y) \ar[rr]_{F(\lambda_Y)}  &&  FT(Y)	&& Y \ar@{ >->}[rr]_{\lambda_Y}  &&  T(Y) \ar@<1ex>[rr]^{\lambda_{T(Y)}} \ar@<-1ex>[rr]_{T(\lambda_{Y})}&& T^2(Y) 
}
$$
then the two right hand side squares of the right hand side diagram commute and, by 1), we get the dotted factorization $k$ in $\EE$ such that $U(h).\lambda_X=\lambda_Y.k$. It remains to show  that $F(k)=h$, which is a consequence of $UF(k).\lambda_X=\lambda_Y.k=U(h).\lambda_X$. The unicity of such a $k$ is a consequence of the fact that $\lambda_{Y}$ is a monomorphism.

Conversely suppose 2). Let $l:Z\to T(X)$ be a map such that $\lambda_{T(X)}.l=T(\lambda_{X}).l$. We get the following commutative square in $\bar\EE$:
$$\xymatrix@=5pt{
	F(Z) \ar[dd]_{\epsilon_{F(X)}.F(l)}  \ar[rr]^{F(\lambda_Z)}  && FT(Z) \ar[dd]^{FU(\epsilon_{F(X)}.F(l))}\\
	&&&&   \\
	F(X) \ar[rr]_{F(\lambda_X)}  &&  FT(X)	 
}
$$
Since: $FU(\epsilon_{F(X)}.F(l)).F(\lambda_{Z})=FU(\epsilon_{F(X)}).FT(l).F(\lambda_{Z})$\\ =$FU(\epsilon_{F(X)}).F(\lambda_{T(X)}).F(l)=F(l)$.\\
While: $F(\lambda_X).\epsilon_{F(X)}.F(l)=\epsilon_{FT(X)}.FT(\lambda_X). F(l)=\epsilon_{FT(X)}.F(\lambda_{T(X)}). F(l)=F(l)$.\\
According to 2), there is a map $k:Z\to X$ such that $F(k)=\epsilon_{F(X)}.F(l)$. Whence $\lambda_X.k=l$, by $\epsilon_{F(X)}.F(\lambda_X.k)=F(k)=\epsilon_{F(X)}.F(l)$. It remains to show the unicity of the factorization. Let $k'$ be be such that $\lambda_X.k'=l$. Since $C_{(U,F)}$ is faithful, checking $k=k'$ is equivalent to checking $F(k)=F(k')$. We get: $F(k')=\epsilon_{F(X)}.F(\lambda_X).F(k')=\epsilon_{F(X)}.F(\lambda_Z.k')=\epsilon_{F(X)}.F(l)=F(k)$.
 The last assertion is straightforward since $F=F_C.C_{(U,F)}$ and $F_C$ is conservative.
\endproof
\begin{prop}\label{comona}
	Let $(U,F): \bar\EE \to \EE$ be  an adjoint pair and $(T,\lambda,\mu)$ its associated monad on $\EE$. Suppose the functor $F$ cartesian. Then the following conditions are equivalent:\\
	1) the natural transformation $\lambda$ is the equalizer of the pair $(\lambda_{T},T(\lambda))$;\\
	2) the functor $F$ is conservative.\\
	Under any of these conditions, the functor $F$ is comonadic.
\end{prop}
\proof
Suppose 1). We shall show that the comparison functor $C_{(U,F)}$ is an equivalence of categories. Accordingly the functor $F$ will be comonadic and thus conservative.
So, let us show that $C_{(U,F)}$ is essentially surjective. Let $a: W\to F.U(W)=C(W)$ be a co-algebra structure on $W$ in $\bar\EE$. By the identities $\epsilon_W.a=1_W$ and $F.U(a).a=F(\lambda_{U(W)}).a$, it produces a $2$-truncated  split simplical object in $\bar\EE$:
$$\xymatrix@=18pt{
	W  \ar@<-2ex>[rr]_{a}  &&	F.U(W) \ar[ll]_{\epsilon_W} 
	\ar@<-4ex>[rrr]_{F.U(a)} \ar[rrr]|{F(\lambda_{U(W)})} &&& (F.U)^2(W)\ar@<-2ex>[lll]_{\epsilon_{F.U(W)}} \ar@<2ex>[lll]|{F.U(\epsilon_W)} 
}
$$
Accordingly, $a$ is the equalizer of the pair $(F.U(a),F(\lambda_{U(W)}))$ in $\bar\EE$. We have to find an object $J$ in $\EE$ such that  $C_{(U,F)}(J)=(F(J),F(\lambda_J)\simeq (W,a)$ in $\ColgC$. For that take the equalizer $j: J\to U(W)$ in $\EE$ of the cosplit parallel pair $(U(a),\lambda_{U(W)})$. Since $F$ is cartesian, this equalizer is preserved by $F$. So the natural comparison $\gamma: F(J)\to W$ such that $a.\gamma=F(j)$ in $\bar\EE$ is an isomorphism. It remains to check that the following square commutes:
$$\xymatrix@=18pt{
	W \ar[r]^{a} & F.U(W) \\
	F(J) \ar[u]^{\gamma} \ar[r]_{F(\lambda_J)} &  F.T(J) \ar[u]_{F.U(\gamma)}
}
$$
namely that $F.U(\gamma).F(\lambda_J)=F(j)$. We shall check it by composition with the monomorphism $F.U(a)$:
$F.U(a).F.U(\gamma).F(\lambda_J)=F.U.F(j).F(\lambda_J)=F(\lambda_{U(W)}).F(j)$ while $F.U(a).F(j)=F(\lambda_{U(W)}).F(j)$ by definition of $j$.

Conversely suppose 2). We have to show that the natural transformation $\lambda$ is the equalizer of the pair $(\lambda_{T},T(\lambda))$. For that take the equalizer $j: J\into X$ in $\EE$ of the cosplit pair $(\lambda_{T(X)},T(\lambda_X))$ and denote $\gamma: X \to J$ the natural comparison such that $\lambda_{X}=j.\gamma$. Since $F$ is cartesian, this equalizer $j$ is preserved by $F$, and, $F\lambda_{X}$ being necessarily the equalizer of the pair $(F(\lambda_{T(X)}),F.T(\lambda_X))$ thanks to the retraction $\epsilon_{F(X)}$, the map $F(\gamma): F(X) \to F(J)$ is an isomorphism is $\bar\EE$. Now, since $F$ is conservative, $\gamma$ is an isomorphism, and   $\lambda_X$ is the equalizer of the pair $(\lambda_{T(X)},T(\lambda_X))$.
\endproof
\begin{coro}\label{auton}
	Let $(T,\lambda,\mu)$ be a monad with a cartesian endofunctor $T$, then the following conditions are equivalent:\\
	1) the monad is half-cartesian;\\
	2) the endofunctor $T$ is conservative.\\
	Under any of these conditions, the adjoint pair $(U^T,F^T): \AlgT\splito \EE$ is an autonomous one.
	
	Moreover, when the monad is hypercartesian, the adjoint pair $(U^T,F^T): \AlgT\splito \EE$ is a cartesian one.
\end{coro}
\proof
The functor $U^T$ being monadic, it reflects pullbacks; so, $T$ is cartesian if and only if $F^T$ is cartesian. Applying the previous proposition, the natural transformation $\lambda$ is the equalizer of the pair $(\lambda_{T},T(\lambda))$ if and only if $F^T$ is conservative. This is the case, $U^T$ being conservative, if and only if $T$ is conservative. Then $F^T$ is comonadic; since $U^T$ is monadic, the adjoint pair $(U^T,F^T): \AlgT\splito \EE$ an autonomous one. The last assertion is a consequence of Proposition \ref{cartadj}.
\endproof

\section{Kleisli category of a monad}

The canonical decomposition of the functor $F^T: \EE \stackrel{\bar F^T}{\longrightarrow} 
KlT \stackrel{K_T}{\longrightarrow} \AlgT$ into a functor $F^T$ which is bijective on objects followed by a fully faithful functor $K_T$ produces the \emph{Kleisli category} $KlT$ \cite{K} of the monad. Accordingly, the functor $\bar F^T$ remains a left adjoint to $:\bar U^T=U^T.K_T: KlT\to \EE$, and obviously this adjoint pair recovers the monad $(T,\lambda,\mu)$ as well. Moreover the functor $\bar U^T$ is conservative as a composition of two conservative functors.

By the adjoint bijection $\Hom_{KlT}(X,Y)\simeq\Hom_\EE(X,TY)$, any map $X \cdots> Y$ in $KlT$ is given by a map $\alpha: X \to T(Y)$ in $\EE$; we call the map $\alpha \in \EE$, the \emph{support} of this map in $KlT$ which we shall then denote by $"\alpha": X \cdots> Y$. Given any other map $"\beta":Y\cdots Z$ in $KlT$, we get $"\beta"."\alpha"="\mu_Z.T(\beta).\alpha"$ in $KlT$. In this way, the natural transformation $\epsilon_X: \bar F^T.\bar U^T(X)=T(X)\cdots> X$ is given by $"1"_{T(X)}$ and, for any map $"\alpha": X\cdots>Y$ in $KlT$, its support $\alpha$ is the unique map in $\EE$ such that
$"\alpha"="1"_{T(X)}.\bar F^T(\alpha)$ in $KlT$. So, given any map $f:X\to Y$ in $\EE$, we get $\bar F^T(f)="\lambda_Y.f"$, and $U^T("\alpha")=\mu_Y.T(\alpha): T(X)\to T(Y)$ in $\EE$. 

We shall now investigate, step by step, how the assumptions of a cartesian monad surprisingly organizes the properties of $KlT$ 

\subsection{Consequences of constraints on $\lambda$}

Let $(T,\lambda,\mu)$ be a monad on $\EE$.

\begin{prop}\label{incl}
The endofunctor $T$ of the monad is faithful if and only if the natural transformation $\lambda$ is monomorphic. The functor $\bar F^T: \EE \to KlT$ is then an inclusion.
\end{prop}
\proof
Suppose $\lambda$ monomorphic. Given a parallel pair $(f,g)$ of maps between $X$ and $Y$ such that $T(f)=T(g)$. Then $T(f).\lambda_X=T(g).\lambda_X$. So, $\lambda_Y.f=\lambda_Y.g$ and $f=g$. 

Conversely suppose the endofunctor $T$ faithful. Given a parallel pair $(f,g)$ of maps between $X$ and $Y$ such that  $\lambda_Y.f=\lambda_Y.g$, we get $T(\lambda_Y).T(f)=T(\lambda_Y).T(g)$.
Since $T(\lambda_Y)$ is a monomorphism as a retract of $\mu_Y$, we get $T(f)=T(g)$; and $f=g$.
The last assertion is then straightforward since  $\bar F T$ is bijective on objects.
\endproof
So, when $\lambda$ is monomorphic, we are in the rather weird situation of a  bijective on objects inclusion $\bar F^T: \EE \into KlT$ which admits a right adjoint $\bar U^T$. The endofunctor of the induced comonad $(C=\bar F^T.\bar U^T,\epsilon,\nu)$ on $KlT$ coincides with $T$ on objects and maps in $\EE$ and we get $C("\alpha")=\mu_Y.T(\alpha)$. Whence the following diagram in $KlT$ where $\mu_X$ coincides with $C("1"_{T(X)})$:
$$\xymatrix@=15pt{
	X  \ar@<-1ex>[rr]_{\lambda_X} \ar@{.>}[dd]_{"\alpha"}  &&	T(X)\ar[dd]_{\mu_Y.\alpha} \ar@{.>}[ll]_{"1"_{T(X)}} 
	\ar@<-2ex>[rr]_{T(\lambda_X)} \ar[rr]|>>>>>{\lambda_{T(X)}} && T^2(X)\ar@{.>}@<-1ex>[ll]_{"1"_{T^2(X)}} \ar@<1ex>[ll]|>>>>>{\mu_X} \ar[dd]^{T(\mu_Y.\alpha)} \\
	&&&&\\ 
	Y  \ar@<-1ex>[rr]_{\lambda_Y}  &&	T(Y) \ar@{.>}[ll]_{"1"_{T(Y)}} 
	\ar@<-2ex>[rr]_{T(\lambda_Y)} \ar[rr]|>>>>>{\lambda_{T(Y)}} && T^2(Y)\ar@{.>}@<-1ex>[ll]_{"1"_{T^2(Y)}} \ar@<1ex>[ll]|>>>>>{\mu_Y}\\
}
$$
Thanks to the horizontal $2$-truncated split simplicial objects, the map $\lambda_X$ appears to be the equalizer in $KlT$ of the pair $(\lambda_{T(X)},T(\lambda_{X}))$.
\begin{prop}\label{wcodj}
	Let $j:\EE\into \bar \EE$ be a bijective on objects inclusion. When $j$ admits a right adjoint $T:\bar\EE \to \EE$, the induced monad $(T,\lambda,\mu)$ on $\EE$ has its $\lambda$ monomorphic. Moreover we get $KlT=\bar\EE$.
\end{prop}
\proof
Let us denote by $\epsilon_X: X<\cdots T(X)$ the co-unit in $\bar\EE$ of this adjunction. The map $1_X:X\to X$ produces a unique map $\lambda_X:X \to T(X)$ in $\EE$
such that $\epsilon_X.\lambda_X=1_X$ and the following diagram in $\bar\EE$:
$$\xymatrix@=18pt{
	X  \ar@<-2ex>[rr]_{\lambda_X}  &&	T(X) \ar@{.>}[ll]_{\epsilon_{X}} 
	\ar@<-4ex>[rr]_{T(\lambda_X)} \ar[rr]|{\lambda_{T(X)}} && T^2(X)\ar@{.>}@<-2ex>[ll]_{\epsilon_{T(X)}} \ar@<2ex>[ll]|{\mu_{X}} 
}
$$
which makes $\epsilon_X$ the coequalizer of the pair $(\epsilon_{T(X)},T(\epsilon_X))$ in $\bar\EE$, and $\lambda_X$ the equalizer of the pair $(\lambda_{T(X)},T(\lambda_X))$ in $\bar\EE$. So, $\lambda_X\in \EE$ is a monomorphism in $\bar\EE$, whence in $\EE$. Since the co-unit $\epsilon_X$ in $\bar \EE$ is the coequalizer of the pair $(\epsilon_{T(X)},T(\epsilon_X))$, the comparison functor $A_{(T,j)}: \bar \EE\to \AlgT$ is fully faithful. Now consider the following commutative diagram:
$$\xymatrix@=18pt{
	\bar\EE  \ar[rr]^{A_{(T,j)}} \ar[dr]_j  &&	\AlgT \ar[dl]_{F^T} \\
	& \EE  \ar@<1ex>[ur]_{U^T}  \ar@<1ex>@{>->}[ul]_T  &
}
$$ 
Since $\EE\into \bar\EE$ is bijective on objects, we get $\bar\EE=KlT$.
\endproof
From Proposition \ref{Tff}, we get the following: 
\begin{coro}\label{Tcons}
	Let $(T,\lambda,\mu)$ be a monad on $\EE$. The two following conditions are equivalent:\\
	1) the natural transformation $\lambda$ is the equalizer of the pair $(\lambda_{T},T(\lambda))$;\\
	2) a map $"\alpha":X\cdots> Y$ in KlT lies in $\EE$ if and only if $\lambda_{T(Y)}.\alpha=T(\lambda_{Y}).\alpha$ in $\EE$.\\
	The inclusion $\bar F^T: \EE\into KlT$ is then conservative; this last point means that  the inverse in $KlT$ of a map $f$ of (the subcategory) $\EE$ belongs to $\EE$.
\end{coro}
\proof
The condition 2) of Proposition \ref{Tff} means that, when the following square commutes in $KlT$, the map $"\alpha"$ belongs to $\EE$:
$$\xymatrix@=5pt{
	X \ar@{.>}[dd]_{"\alpha"}  \ar[rr]^{\lambda_X}  && T(X) \ar[dd]^{\mu_Y.T(\alpha)}\\
	&&&&   \\
Y \ar[rr]_{\lambda_Y}  &&  T(Y)	 
}
$$ 
The commutation in $KlT$ is: $T(\lambda_{Y}).\alpha=\lambda_{T(Y)}.\mu_Y.T(\alpha).\lambda_{X}$ in $\EE$; now this last term is clearly: $\lambda_{T(Y)}.\mu_Y.\lambda_{T(Y)}.\alpha=\lambda_{T(Y)}.\alpha$. When 1) is true, the endofunctor $T$ is conservative; since $T=\bar U^T.\bar F^T$, the functor $\bar F^T$ is conservative as well.
\endproof
Let us introduce the following:
\begin{defi}
	A class $\Sigma$ of maps in a category $\EE$ is said to be left cancellable, when it is such that $g.f\in \Sigma$ and $g\in \Sigma$ imply $f\in\Sigma$.
\end{defi}
\begin{prop}\label{c3from2}
	Let $(T,\lambda,\mu)$ be a monad with $\lambda$ monomorphic. Then the following conditions are equivalent:\\
	1) the natural transformation $\lambda$ is cartesian;\\
	2) the inclusion $\bar F^T: \EE \into KlT$ makes $\EE$ a left cancellable subcategory of $KlT$.\\
	Any of these conditions implies that $\lambda$ is the equalizer of the pair $(\lambda_{T},T(\lambda))$. 
\end{prop}
\proof
Suppose 1). Let $h:X \to Z \in \EE, \; g:Y\to Z \in\EE$, and $"\phi": X\cdots>Y \in KlT$ be such that $g."\phi"=h$ in $KlT$. This means that $T(g).\phi=\lambda_Z.h$ in $\EE$. According to the previous corollary, we must show that $\lambda_{T(Y)}.\phi=T(\lambda_{Y}).\phi$ in $\EE$. Now consider the following diagram in $\EE$:
$$\xymatrix@=8pt{
	X  \ar[rr]^{\phi} \ar[dd]_{h}  &&	T(Y)\ar[dd]_{T(g)}  \ar@<-1ex>[rr]_{T(\lambda_Y)}  
	\ar@<1ex>[rr]^{\lambda_{T(Y)}}  && T^2(Y)  \ar[dd]^{T^2(g)} \ar[ll]|{\mu_Y} \\
	&&&&\\ 
	Z  \ar[rr]_{\lambda_Z}  &&	T(Z) \ar@<1ex>[rr]^{\lambda_{T(Z)}} 
	\ar@<-1ex>[rr]_{T(\lambda_Z)}  && T^2(Z)\ar[ll]|{\mu_Z} \\
}
$$
Since the leftward right hand side commutative square is a pullback, it is enough to check our equality via composition with $\mu_Y$ (trivial) and by $T^2(g)$ which is straightforward since the left hand side square commutes.
 
Conversely suppose 2). We have to show that the following left hand side square is a pullback in $\EE$. So let $(h,\phi)$ be a pair of maps in $\EE$ such that $T(g).\phi=\lambda_Z.h$. This means that the following right hand side triangle commutes in $KlT$:
$$
\xymatrix@=5pt{
 Y  \ar[rr]^{\lambda_Y} \ar[dd]_{g}  &&	T(Y)\ar[dd]^{T(g)}	&&&&&& Y \ar[ddrr]^{g}\\
	&&&&&&   \\
 Z  \ar[rr]_{\lambda_Z}  &&	T(Z)	&&&& X \ar[rrrr]_{h} \ar@{.>}[rruu]^{"\phi"} &&&&  Z
}
$$
According to our assumption, the map $"\phi"$ is in $\EE$, which means that there is a map $f: X\to Y$ in $\EE$ such that $\phi=\lambda_Y.f$ and $g.f=h$ in $\EE$. This map $f$ is unique since $\lambda_Y$ is monomorphic. So, the square in question is a pullback.
\endproof
Finally we get:
\begin{prop}\label{lcwcodj}
	Let $j:\EE\into \bar \EE$ be a bijective on objects left cancellable inclusion. When $j$ admits a right adjoint $T:\bar\EE \to \EE$, the induced monad $(T,\lambda,\mu)$ on $\EE$ is such that  $\lambda$ is cartesian.
\end{prop}
\proof
According to Proposition \ref{wcodj}, $\lambda_X$ is monomorphic and $\bar\EE=KlT$.
Since $\EE$ is left cancellable in $\bar\EE$, the previous proposition asserts that $\lambda$ is cartesian.
\endproof

\subsection{Consequences of constraints on $T$}

Rephrasing Corollary \ref{auton} we get:
 
\begin{coro}\label{onT}
	Let $(T,\lambda,\mu)$ be a monad with a cartesian endofunctor $T$. When $\lambda$ is the equalizer of the pair $(\lambda_T,T\lambda)$, the injective functor $\bar F^T: \EE \into KlT$ is conservative, cartesian and then it reflects the pullbacks of $\EE$. 
\end{coro}

\subsection{Consequences of constraints on $\mu$}

\begin{defi}
	A pullback stable class $\Sigma$ of maps in a category $\EE$ is a class of maps which admit pullbacks along any map in $\EE$ and whose pullbacks stay in $\Sigma$.
\end{defi}

\begin{prop}\label{spullback}
	Let  $(T,\lambda,\mu)$ be a monad with a cartesian endofunctor $T$ and a cartesian natural transformation $\mu$.	Then, the class $\bar F^T(\Sigma)$ is pullback stable in $KlT$.
\end{prop}
\proof
Starting with a map $f$, consider the following left hand side pullback in $\EE$:
$$
\xymatrix@=7pt{
	 V \ar[dd]_{\phi}  \ar[rr]^{h}  && U \ar[dd]^{\psi} && V \ar@{.>}[dd]_{"\phi"}  \ar@{.>}[rrr]^{\bar F^T(h)}  &&& U \ar@{.>}[dd]^{"\psi"}\\
	&&&&   \\
 T(X) \ar[rr]_{T(f)}  &&  T(Y) && X \ar@{.>}[rrr]_{\bar F^T(f)}  &&&  Y
}
$$
This pullback in $\EE$ determines a right hand side commutative square in $KlT$. Let us show that it is a pullback in this category. Since $K_T: KlT \to \AlgT$ is fully faithful it is sufficient to check it in $\AlgT$; now, since $U^T$ reflects pullbacks, it is sufficient to check that its image by $U^T$, which is given the following vertical rectangle, is a pullback in $\EE$:
$$
\xymatrix@=10pt{
	T(V) \ar[d]_{T(\phi)}  \ar[rr]^{T(h)}  && T(U) \ar[d]^{T(\psi)} \\
	T^2(X) \ar[d]_{\mu_X} \ar[rr]^{T^2(f)} &&  T^2(Y)\ar[d]^{\mu_Y} \\
	T(X) \ar[rr]_{T(f)}  &&  T(Y) 
}
$$
This is the case since the lower square is a pullback, $\mu$ being cartesian, and the upper one as well since, $T$ being cartesian, the image by $T$ of our above left hand side pullback is preserved by $T$. 
\endproof
From Propositions \ref{c3from2}, \ref{spullback} and \ref{comona}, we get:
\begin{coro}\label{finalKlT}
	Let $(T,\lambda,\mu)$ be a cartesian monad on $\EE$. We then get a bijective on objects inclusion $\bar F^T: \EE \into KlT$, which makes $\EE$ a left cancellable subcategory of $KlT$ which is pullback stable in $KlT$. This inclusion functor is cartesian and conservative. It admits a right adjoint $\bar U^T$ which necessarily makes this inclusion comonadic. 
\end{coro}
And from Proposition \ref{lcwcodj} the following:
\begin{coro}
	Let $j:\EE\into \bar\EE$ be a bijective on objects left cancellable and pullback stable inclusion. When $j$ admits a right adjoint $T:\bar\EE \to \EE$, the induced monad $(T,\lambda,\mu)$ on $\EE$ is cartesian and $\bar\EE$ is the Kleisli category of this monad.
\end{coro}
\proof
By Proposition \ref{wcodj}, we know that $\lambda$ is cartesian and $\bar\EE=KlT$. The endofunctor $T.j:\EE\to\EE$ is cartesian since so are $T$ (being a right adjoint) and $j$ ($\EE$ being pullback stable in $KlT$, the inclusion $j$ preserves the pullbacks). It remains to check that $\mu$ is cartesian. For that, given any map $f\in\EE$, consider the following leftward pullback square in $KlT$:
$$
\xymatrix@=7pt{
	&&&& T(X) \ar@{.>}@<-1ex>[lllld]_{"1_{T(X)}"} \ar[llddd]^{T(f)} \ar[dll]^>>>>{\check f}\\
	X \ar[dd]^{f}   && P \ar[dd]^{\bar f} \ar@{.>}[ll]^{"\phi_f"}\\
	&&&&&&   \\
	Y   &&  T(Y) \ar@{.>}[ll]^{"1_{T(Y)}"}
}
$$
The commutation of this square in $KlT$ means $T(f).\phi_f=\bar f$ in $\EE$. The commutation of the quadrangle in $KlT$ produces a factorization $\check f$ in $KlT$; since $\EE$ is left cancellable in $KlT$, the commutation of the vertical right hand side triangle makes $\check f$ in $\EE$. So we get $\bar f.\check f=T(f)$ in $\EE$ and $"\phi_f".\check f="1_{T(X)}"$ in $KlT$ which means  $\phi_f.\check f=1_{T(X)}$ in $\EE$. Checking $\phi_f.\check f=1_{T(X)}$ in $\EE$ will prove that the quadrangle is a pullback in $KlT$, which, being preserved by $T$, will prove, in turn, that $\mu$ is cartesian. We shall check $\phi_f.\check f=1_{T(X)}$, by composition with $\bar f$ and $"\phi_f"$ in $KlT$. 1): $\bar f.\phi_f.\check f=T(f).\check f=\bar f=\bar f.1_{T(X)}$; 2)$"\phi_f".\phi_f.\check f="\phi_f".1_{T(X)}$ in $KlT$ is equivalent to $\phi_f.\check f.\phi_f=\phi_f$ in $\EE$, which is true since $\phi_f.\check f=1_{T(X)}$.
\endproof
 
\section{Internal categories}\label{ccategory}

\emph{From now on, we shall suppose any ground category $\EE$ with pullbacks and terminal object $1$}. Given any map $f$, we shall denote the kernel equivalence $R[f]$ of this map (which is underlying an internal groupoid $R[f]_{\bullet}$ in $\EE$) in the following way:
$$\xymatrix@=20pt
{
	R_2[f]  \ar@<2ex>[rr]^{p_0^f} \ar@<-2ex>[rr]_{p_2^f} \ar[rr]|{p_1^f}&&	R[f]  \ar@<2ex>[rr]^{p_0^f} \ar@<-2ex>[rr]_{p_1^f} && X  \ar[rr]^{f} \ar[ll]|{s_0^f} && Y
}
$$ 
and given any commutative square, as on the right hand side,
we denote by $R(\phi)$ the induced map between the kernel equivalences:
$$\xymatrix@=20pt
{
	R[f] \ar[d]_-{R(\phi)} \ar@<1ex>[r]^{p_0^f} \ar@<-1ex>[r]_{p_1^f} & X  \ar[d]^-{\phi} \ar[r]^{f} \ar[l]|{s_0^f} &
	Y \ar[d]^-{\psi} \\
	R[f'] \ar@<1ex>[r] \ar@<-1ex>[r] & X' \ar[r]_-{f'} \ar[l] & Y'.
}
$$
which is underlying an internal functor $R(\phi)_{\bullet}:R[f]_{\bullet}\to R[f']_{\bullet}$. Let us recall the following useful Barr-Kock Observation \cite{BG}:
\begin{lemma}\label{BK}
	When the above right hand side square is a pullback, the left hand part of the diagram is a discrete fibration between groupoids, which implies that any vertical commutative square is a pullback. Conversely, if the left hand side part of the diagram is a discrete fibration, and if $f$ is a pullback stable strong epimorphisms (it is the case when $f$ is split), then the right hand side square is a pullback.
\end{lemma}

\subsection{Basics}

Internal categories have been introduced by Ch. Ehrhesmann in \cite{Eh0}. We deliberately choose the simplicial notations. For the basics on simplicial objects, see, for instance, Chapter VII.5 in \cite{ML}. An internal category in $\EE$ is is a $3$-truncated simplicial object $X_{\bullet}$ in $\EE$, namely a diagram:
$$\xymatrix@=20pt
{
	X_{\bullet}:&& X_3 \ar@<6ex>[rr]|>>>>>>>>{d_{3}^{X_{3}}} 
	\ar@<2ex>[rr]|>>>>>>>>{d_{2}^{X_{3}}} \ar@<-2ex>[rr]|>>>>>>>>{d_{1}^{X_{3}}}\ar@<-6ex>[rr]_>>>>>>>>{d_{0}^{X_{3}}}	&& X_2  \ar@<4ex>[ll]|<<<<<<{s_0^{X_{3}}} \ar[ll]|<<<<<<{s_1^{X_{3}}} \ar@<-4ex>[ll]|<<<<<<{s_2^{X_{3}}} \ar@<4ex>[rr]^{d_{2}^{X_{2}}} \ar[rr]|{d_{1}^{X_{2}}} \ar@<-4ex>[rr]_{d_{0}^{X_{2}}} && X_1 \ar@<2ex>[rr]^{d_{1}^{X_{1}}} \ar@<-2ex>[rr]_{d_{0}^{X_{1}}} \ar@<-2ex>[ll]|<<<<<<{s_1^{X_{2}}} \ar@<2ex>[ll]|<<<<<<{s_0^{X_{2}}}&&
	X_0   \ar[ll]|{s_0^{X_{1}}}  
}
$$
(where we shall drop the upper indexes when there is no ambiguities)
 subject to the following identities:
$$
\begin{aligned}
d_i.d_{j+1}&=d_j.d_i, \;\;\; i\leq j  & \;\;\; d_i.s_j&=s_{j-1}.d_i, \;\;\; i<j \\
s_{j+1}.s_i&=s_i.s_j, \;\;\; i\leq j  & \;\;\; d_i.s_j&=1, \;\;\;\;\;\;\;\;\;\;\;\; i=j, j+1\\ 
&                           & \;\;\; d_i.s_j&=s_j.d_{i-1}, \;\;\; i>j+1  
\end{aligned}
$$
where the object $X_2$ (resp. $X_3$) is obtained by the pullback of $d_0^{X_1}$ along $d_1^{X_1}$ (resp. $d_0^{X_2}$ along $d_2^{X_2}$).  An internal functor is a simplicial morphism between this kind of $3$-truncated simplical objects. Let us recall some classical classes of internal functors:
\begin{defi}\label{dfib}
An internal functor $f_{\bullet}: X_{\bullet}\to Y_{\bullet}$ is a discrete cofibration (resp. fibration) when the following right hand side square horizontally indexed by $0$ (resp. by $1$) is a pullback:
$$\xymatrix@=12pt{
	 X_2 \ar[dd]_{f_2}\ar@<-2ex>[rr]_{d_2} \ar[rr]|{d_1} \ar@<2ex>[rr]^{d_0}  &&	X_1 \ar[dd]^{f_1}\ar@<-2ex>[rr]_{d_1} 
	\ar@<2ex>[rr]^{d_0}  && X_0\ar[ll]|{s_0} \ar[dd]^{f_0}\\
	&&&&\\
   Y_2\ar@<-2ex>[rr]_{d_2} \ar[rr]|{d_1} \ar@<2ex>[rr]^{d_0}  &&	Y_1 \ar@<-2ex>[rr]_{d_1} 
	\ar@<2ex>[rr]^{d_0}  && Y_0\ar[ll]|{s_0}
}
$$
\end{defi} 
We denote by $Cat\EE$ the category of internal categories in $\EE$, and by $(\;)_0: Cat\EE\to \EE$ the 
forgetful functor associating with any internal category $X_{\bullet}$ its ''object of objects" $X_0$. Since $\EE$ has pullbacks, so is the category $Cat\EE$ since, by commutation of limits, it is easy to see that the limits in $Cat\EE$ are built levelwise in $\EE$. So, the forgetful functor $(\;)_0$ is cartesian.

The functor $(\;)_0$ is actually a fibration whose cartesian maps are the internal \emph{fully faithful functors} (obtained by a joint pullback) and whose maps in the fibers are the internal functors which are "identities on objects" (\emph{ido-functors} or \emph{idomorphisms} for short).

It is clear that the fiber $Cat_1\EE$ above the terminal object $1$ is nothing but the pointed category $Mon\EE$ of internal monoids in $\EE$. Any fiber $Cat_Y\EE$ above an  object $Y$, with  $Y\neq 1$, has an initial object with the discrete equivalence relation $\Delta_Y=R[1_Y]$ and a terminal one with the indiscrete one $\nabla_Y=R[\tau_Y]$, where $\tau_Y: Y\to 1$ is the terminal map. So, the left exact fully faithful functor $\nabla: \EE \to Cat\EE$ admits the fibration $(\;)_0$ as left adjoint and makes the pair $((\;)_0,\nabla)$ a \emph{fibered reflection} in the sense of \cite{B-1} (see also section \ref{locfib} below). A functor $f_{\bullet}$ is then cartesian with respect to $(\;)_0$ (namely, internally fully faithful) if and only if the following left hand side square is a pullback in $Cat\EE$, or, equivalently  the right hand side one is a pullback in $\EE$:
$$\xymatrix@=20pt{
	X_{\bullet}  \ar[r]^{f_{\bullet}} \ar[d]_{} & Y_{\bullet}  \ar[d]^{} && X_1 \ar[r]^{f_1} \ar[d]_{(d_0,d_1)} & Y_1  \ar[d]^{(d_0,d_1)}\\
	\nabla_{X_{\bullet}}   \ar[r]_{\nabla_{f_{\bullet}}}  & \nabla_{Y_{\bullet}}  &&  X_0\times X_0   \ar[r]_{f_0\times f_0}  & Y_0\times Y_0 }
$$
As for any left exact fibration, we get:
\begin{prop}\label{prop1}
	1) The cartesian maps (= internally fully faithful) functors are stable  under composition and pullback.
	
	2) Given any commutative square in $Cat\EE$ where both $x_{\bullet}$ and $y_{\bullet}$ are cartesian maps:
	$$\xymatrix@=20pt{
		X_{\bullet}  \ar[r]^{x_{\bullet}} \ar@<-2pt>[d]_{f_{\bullet}} & X'_{\bullet}  \ar@<-1pt>[d]^{f'_{\bullet}} \\
		Y_{\bullet}   \ar[r]_{y{\bullet}}  & Y'_{\bullet} }
	$$
	then it is a pullback:\\
	1) if and only if its image by $(\;)_0$ is a pullback\\
	2) in particular when $f_{\bullet}$ and $f'_{\bullet}$ are ido-functors.
\end{prop}
The dual category $X^{op}_{\bullet}$ of $X_{\bullet}$ is the internal category where the role of $d_0^{X_1}$ and $d_1^{X_1}$ are interchanged. In the context of monads, we have the following:

\begin{prop}\label{Mum}
	Let $(T,\lambda,\mu)$ be any monad on $\EE$ with $\mu$ cartesian. Then, for any algebra $(X,\xi)$, the following diagram produces an internal category 
 $\overline T(X,\xi)$  in $\AlgT$:
	$$\xymatrix@=18pt{
		(T^3(X),\mu_{T^2(X)})\ar@<-4ex>[rr]_{T^2(\xi)} \ar[rr]|{T(\mu_X)} \ar@<4ex>[rr]^{\mu_{T(X)}}  &&	(T^2(X),\mu_{T(X)}) \ar@<-2ex>[rr]_{T(\xi)} 
		\ar@<2ex>[rr]^{\mu_X} \ar@<2ex>[ll]|{T^2(\lambda_X)} \ar@<-2ex>[ll]|{T(\lambda_{T(X)})} && (T(X),\mu_X )\ar[ll]|{T(\lambda_X)}
	}
	$$
determining the following commutative square:
$$\xymatrix@=5pt{
	 \Alg T \ar[dd]_{U^T} \ar[rr]^{\bar T}  && Cat\Alg T \ar[dd]^{(\;)_0}\\
	&&	&&   \\
 \EE \ar[rr]_{F^T}  && \Alg T 
}
$$
Moreover any morphism $f: (X,\xi)\to (Y,\gamma)$ of $T$-algebras produces a discrete fibration $\bar T(f): \bar T(X,\xi)\to \bar T(Y,\gamma)$ in $Cat\Alg T$.
\end{prop}
\proof
First, we have the following pullbacks in $\EE$ since $\mu$ is cartesian:
$$\xymatrix@=7pt{
T^4(X) \ar[dd]_{\mu_{T^2(X)}} \ar[rr]^{T^3(\xi)} && T^3(X) \ar[dd]_{\mu_{T(X)}} \ar[rr]^{T^2(\xi)}  && T^2(X) \ar[dd]^{\mu_{X}}\\
&&	&&   \\
T^3(X) \ar[rr]_{T^2(\xi)} && T^2(X) \ar[rr]_{T(\xi)}  && T(X) 
}
$$
Then $T(\mu_X)$, thanks to the axioms of $T$-algebra on $X$, furnishes a composition map. The following pullback in $\EE$:
$$\xymatrix@=5pt{
	 T^2(X) \ar[rr]^{\mu_{X}} \ar[dd]_{T^2(f)}  && T(X) \ar[dd]^{T(f)}\\
	&&   \\
 T^2(Y) \ar[rr]_{\mu_{Y}}  && T(Y) 
}
$$
determines the last assertion.
\endproof
There is a comonad $(Dec,\epsilon,\nu)$ on the simplicial objects which is stable on $Cat\EE$ as soon as $\EE$ has pullbacks. We shall briefly describe this endofunctor $Dec$ and the co-unit $\epsilon: Dec \to 1_{Cat\EE}$ because they will be useful later on. Let us start with the lower internal category $X_{\bullet}$ and  consider the following diagram, where $X_4$ is determined by the pullback of $d_3:X_3\to X_2$ along $d_0: X_3\to X_2$:
$$\xymatrix@=15pt{
DecX_{\bullet}\ar[dd]_{\epsilon_{X_{\bullet}}}: & X_4\ar[dd]_{d_4} \ar@<-3ex>[rr]_{d_3} \ar@<-1ex>[rr]|{d_2} \ar@<1ex>[rr]|{d_1} \ar@<3ex>[rr]^{d_0}	&& X_3 \ar[dd]_{d_3}\ar@<-2ex>[rr]_{d_2} \ar[rr]|{d_1} \ar@<2ex>[rr]^{d_0}  &&	X_2 \ar[dd]^{d_2}\ar@<-2ex>[rr]_{d_1} 
	\ar@<2ex>[rr]^{d_0}  && X_1\ar[ll]|{s_0} \ar[dd]^{d_1}\\
	&&&&\\
X_{\bullet}: & X_3 \ar@<-3ex>[rr]_{d_3} \ar@<-1ex>[rr]|{d_2} \ar@<1ex>[rr]|{d_1} \ar@<3ex>[rr]^{d_0}	&&	X_2\ar@<-2ex>[rr]_{d_2} \ar[rr]|{d_1} \ar@<2ex>[rr]^{d_0}  &&	X_1 \ar@<-2ex>[rr]_{d_1} 
	\ar@<2ex>[rr]^{d_0}  && X_0\ar[ll]|{s_0}
}
$$
The category $DecX_{\bullet}$ is given by the upper row (in $Set$, it is the sum of all the coslice categories $Y/\EE$), while the co-unit of the comonad is given by the vertical internal functor which is a discrete cofibration.

\subsection{The cartesian monad $(T_{X_{\bullet}},\lambda_{X_{\bullet}},\mu_{X_{\bullet}})$ on $\EE/X_0$}\label{T_X}

Given any objet $Y$, the slice category $\EE/Y$ has the maps with codomain $Y$ as objects, and the commutative triangles above $Y$ as morphisms. Given any map $g:Z\to Y$, the composition with $g$ determines a functor $\Sigma_g:\EE/Z\to \EE/Y$ which admits as right adjoint the pullback functor $g^*: \EE/Y\to \EE/Z$ along $g$.

Let $X_{\bullet}$ be an internal category in $\EE$. It produces a cartesian monad on the slice category $\EE/X_0$ whose endofunctor is $T_{X_{\bullet}}=\Sigma_{d_0}.d_1^*$. The following diagram where any leftward plain square is a pullback describes vertically the behaviour of the endofunctor $T_{X_{\bullet}}$  from left to right. The associated natural transformations $\lambda_{X_{\bullet}}$ and $\mu_{X_{\bullet}}$ are precisely described by the upper horizontal dotted arrows $\sigma_0^h$ and $\delta_1^h$  which are induced by the middle horizontal ones:
$$\xymatrix@=20pt{
	Z  \ar@{ >.>}[r]^{\sigma_0^h} \ar[d]_{h} & d_1^*(Z) \ar@<1ex>[l]^{\delta_1^h} \ar@<-2pt>[d]^{d_1^*(h)} & (d_1.d_2)^*(Z) \ar@<-3pt>[d]^{(d_1.d_2)^*(h)}\ar@<1ex>[l]^{\delta_2^h} \ar@<-1ex>@{.>}[l]_{\delta_1^h}\\
	X_0   \ar@{ >.>}[r]^{s_0}  & X_1 \ar[d]^{d_0}\ar@<1ex>[l]^{d_1} & X_2 \ar[d]^{d_0}\ar@<1ex>[l]^{d_2}\ar@<-1ex>@{.>}[l]_{d_1}\\
	 & X_0 & X_1 \ar@<1ex>[l]^{d_1}\ar[d]^{d_0}\ar@<-1ex>@{.>}[l]_{d_0} \\
	 & & X_0 \\
	h \ar@{ >->}[r]_{\sigma_0^h} & {\;\;T_{X_{\bullet}}(h)} & {\;\;\;T^2_{X_{\bullet}}(h)} \ar[l]^{\delta_1^h}
}
$$
It is a cartesian monad on $\EE/X_0$ since the upper plain part of the diagram is made of pullbacks. It is well known (see for instance \cite{Jon}) that the algebras of this monad coincides with the internal discrete fibrations above $X_{\bullet}$; namely, we have $\Alg T_{X_{\bullet}}=DisF/X_{\bullet}$.

\section{Internal groupoids}\label{groupoid}

\subsection{Basics}

A category $X_{\bullet}$ is a groupoid if and only if any map is invertible. It is a property, which, internally speaking, is equivalent to saying that the following square in the $3$-truncated simplicial object defining $X_{\bullet}$ is a pullback:
	$$\xymatrix@=20pt{
	X_2  \ar[d]_{d_0} \ar[r]^{d_1} & X_1  \ar[d]^{d_0} \\
	X_ 1  \ar[r]_{d_0}  & X_0 
}
$$
It is easy to check (via the Yoneda embedding) that a category $X_{\bullet}$ is a groupoid if and only if any commutative square in the $3$-truncated simplicial object of its definition is a pullback, see \cite{B-1}. The category $Grd\EE$ of internal groupoids is the full subcategory of $Cat\EE$ whose objects are the groupoids, and it determines a sub-fibration:
	$$\xymatrix@=20pt{
	Grd\EE  \ar@{ >->}[r] \ar@<-2pt>[d]_{(\;)_0} & Cat\EE  \ar@<-2pt>[d]^{(\;)_0} \\
\EE   \ar@{=}[r]  & \EE
}
$$
The fibre $Grd_Y\EE$ has the same initial and terminal objects as the fibre $Cat_Y\EE$.

\subsection{The fibration of points}

We denote by $\Pt\EE$  the category whose objects are the split epimorphisms $(g,t):X\splito Y$ in $\EE$ and whose morphisms are the commutative squares between them:
$$\xymatrix@=8pt
{
	X  \ar[dd]_{g} \ar[rr]^{x} && X' \ar[dd]_{g'}	\\
	&&& \\
	Y \ar[rr]_{y} \ar@<-1ex>[uu]_{t} 	&& Y'  \ar@<-1ex>[uu]_{t'} 
}
$$
We denote by $\P_{\EE}: \Pt\EE\to\EE$ the functor which associates with any split epimorphism $(g,t)$ its codomain $Y$, and associates with any morphism $(y,x)$ the map $y$. It is a left exact fibration whose cartesian maps are those ones such that the above square is a pullback of split epimorphisms in $\EE$; it is called the \emph{fibration of points} \cite{B0}. The class $\P$ of cartesian maps is stable under composition and pullbacks in $\Pt\EE$, it is left cancellable and contains the isomorphisms. Accordingly, it determines a bijective on objects inclusion $j_\P: \Sigma_\P\into \Pt\EE$, where $\Sigma_\P$ denotes the subcategory of $\Pt\EE$ whose morphisms belong to the class $\P$, which is left cancellable and pullback stable in $\Pt\EE$.

The fibre above $Y$ is denoted by $\Pt_Y\EE$ and an object of this fibre is called a (generalized) \emph{point} of $Y$, while any morphism in a  fiber is, for short, called an \emph{idomorphism} (=having an identity as lower map $y$). The left exact change of base functor produced by the map $\psi:Y\rightarrow Y'\in\EE$ is the pullback along it and denoted by:
$\psi^{*}:\Pt_{Y'}\mathbb E \rightarrow \Pt_Y\mathbb E$.

\subsection{The monad $(\GG,\sigma,\pi)$ on $\Pt\EE$}

The endofunctor $\GG$ on $\Pt\EE$ defined by $\GG(g,t)=(p_0^g,s_0^g)$ is underlying a monad described by the following diagram in $\mathbb E$, where $t_1=(t.g,1_X)$:
$$\xymatrix@=6pt
{
 X  \ar[dd]_{g} \ar@{.>}[rr]^{t_1} && R[g] \ar[dd]_{p_0^g}	&& R_2[g] \ar@{.>}[ll]_{p_2^g} \ar[dd]_{p_0^g} \\
 &&& \\
 Y \ar@{.>}[rr]_{t} \ar@<-1ex>[uu]_{t} && X   \ar@<-1ex>[uu]_{s_0^g} 
 	&& R[g]  \ar@<-1ex>[uu]_{s_0^g} \ar@{.>}[ll]^{p_1^g}\\
 	(g,t) \ar@{.>}[rr]_{\sigma_{(g,t)}} && \GG(g,t) && \ar@{.>}[ll]^{\pi_{(g,t)}}\GG^2(g,t)
}
$$
It is clear that the maps $\sigma_{(g,t)}$ and $\pi_{(g,t)}$ belong to the class $\P$ of cartesian maps. Although being not strictly cartesian, this monad shares many properties with this notion. Given any pullback stable class $\Sigma$ of maps in a category $\EE$, we shall call \emph{$\Sigma$-cartesian} any functor or natural transformation which only satisfies the cartesian condition on the maps in $\Sigma$:
\begin{defi}
	Given any pullback stable class $\Sigma$ of maps in a category $\EE$ and any monad $(T,\lambda,\mu)$ on $\EE$, we shall say that:\\
	1) this monad is $\Sigma$-cartesian, when:\\
	i) the endofunctor $T$ preserves the maps in $\Sigma$ and is $\Sigma$-cartesian ;\\
	ii) the natural transformations $\lambda$ is the equalizer of the pair $(\lambda_{T},T(\lambda))$;\\
	iii) the natural transformations $\mu$ is $\Sigma$-cartesian;\\
	2) this monad is  strongly $\Sigma$-cartesian when:\\
	i)  it is $\Sigma$-cartesian;\\
	ii)  any $\lambda_X$ and any $\mu_X$ belong to $\Sigma$.
\end{defi}
\noindent\textbf{Remark.} \emph{1) A $\Sigma$-cartesian monad is such that $\lambda$ is $\Sigma$-cartesian;\\2) a monad is a strongly $\Sigma$-cartesian one if and only if:\\
i) the endofunctor $T$ preserves the maps in $\Sigma$ and is $\Sigma$-cartesian ;\\
ii) the natural transformations $\lambda$ and $\mu$ are $\Sigma$-cartesian;\\
iii)   any $\lambda_X$ and any $\mu_X$ belongs to $\Sigma$.}.
\proof
The point 1) is obtained by the proof of Proposition \ref{lcart} restricted to maps in $\Sigma$.  As soon as $\lambda_X$ belongs to $\Sigma$ and $\lambda$ is $\Sigma$-cartesian, then  the natural transformations $\lambda$ is the equalizer of the pair $(\lambda_{T},T(\lambda))$, whence 2).
\endproof

\begin{prop}\cite{B-1}\label{rich}
	The endofunctor $\GG$ is cartesian. It preserves and reflects the maps of the class $\P$. The monad $(\GG,\sigma,\pi)$ is strongly $\P$-cartesian. Furthermore, given any object $(g,t)$ in $\Pt\mathbb E$, the following diagram is a kernel equivalence relation in $\Pt\EE$ with its (levelwise) quotient:
	$$\xymatrix@=25pt
	{
		\GG^3(g,t) \ar@<-2ex>[r]_{\GG\pi_{(g,t)}} \ar@<3ex>[r]^{\pi_{\GG(g,t)}} & \GG^2(g,t) \ar[r]_{\pi_{(g,t)}} \ar[l]_{\GG\sigma_{\GG(g,t)}}  & \GG(g,t) 
	}
	$$
	\end{prop}
\proof The endofunctor $\GG$ is cartesian because it is the result of a cartesian construction. The second assertion comes from Lemma \ref{BK}.
 We already noticed that the maps $\sigma_{(g,t)}$ and $\pi_{(g,t)}$ belong to the class $\P$.
It remains to show that the natural transformations $\sigma$ and $\pi$ are $\P$-cartesian. The second point is a consequence of the fact that the pullbacks in $\Pt\EE$ are levelwise  and of Lemma \ref{BK}, while the first one is a consequence of the further fact that, in the pullback of a split epimorphism, the splitting is pullbacked as well. The last assertion comes from the fact that in the diagram of a kernel equivalence relation:
$$\xymatrix@=15pt
{
	R_2[f]  \ar@<2ex>@{.>}[rr]^{p_0^f} \ar@<-2ex>[rr]_{p_2^f} \ar[rr]|{p_1^f }&&	R[f]  \ar@<2ex>@{.>}[rr]^{p_0^f} \ar@<-2ex>[rr]_{p_1^f} && X  \ar@{.>}[rr]^{f} \ar@{.>}[ll]|{s_0^f} && Y
}
$$
the diagram in plain arrow is again a kernel equivalence relation.
\endproof
Again from \cite{B-1}, recall the following:
\begin{theo}\label{carac}
	Any algebra $\alpha:\GG(g,t)\rightarrow (g,t)$ of this monad is necessarily a $\P$-cartesian map.
	The category of algebras of the monad $(\GG,\sigma,\pi)$ on $\Pt\mathbb E$ is the category $Grd\mathbb E$ of internal groupoids in $\mathbb E$. The functor $U^{\GG}: Grd\EE \to \EE$ associates $(d_0,s_0):X_1\splito X_0$ with any groupoid $X_{\bullet}$, while the functor $F^{\GG}: \EE \to Grd\EE$ associates the internal groupoid $R_{\bullet}[g]$ with split epimorphism $(g,t)$.
\end{theo}
\proof
We refer to \cite{B-1} for the proof. The main conceptual aspects of the proof are detailed in the proof of Proposition \ref{ccartalg}.
\endproof
Accordingly an internal groupoid is given by a reflexive graph in $\EE$ as on the right hand side, completed by a map $d_2:R[d_0]\to X_1$ which makes the following diagram :
$$\xymatrix@=15pt{
	X_{\bullet}: & R_2[d_0] \ar@<-3ex>@{.>}[rr]_{d_3=R[d_2]} \ar@<-1ex>[rr]|{p_2} \ar@<1ex>[rr]|{p_1} \ar@<3ex>[rr]^{p_0}	&&	R[d_0]\ar@<-2ex>@{.>}[rr]_{d_2} \ar[rr]|{p_1} \ar@<2ex>[rr]^{p_0}  &&	X_1 \ar@<-2ex>[rr]_{d_1} 
	\ar@<2ex>[rr]^{d_0} \ar@<1ex>[ll]|{s_1} \ar@<-1ex>[ll]|{s_0}&& X_0\ar[ll]|{s_0}
}
$$
a $3$-truncated simplicial object; namely the map $d_2$ must: 1) complete a $2$-simplicial object, and 2) satisfy $d_2.R(d_2)=d_2.p_2$.
In the set theoretical context, the map $d_2$ is defined by $d_2(\alpha,\beta)=\beta.\alpha^{-1}$ for any pair $(\alpha,\beta)$ of arrows with same domain.
\endproof
\begin{coro}\label{GGaut}
The adjoint pair $(U^\GG,F^\GG): Grd\EE \splito \EE$ is autonomous and the comonad induced on $Grd\EE$ coincides with the rectriction to $Grd\EE$ of the comonad $(Dec,\epsilon,\nu)$ on $Cat\EE$. The endofunctor $Dec: Cat\EE\to Cat\EE$ reflects the groupoids. 
\end{coro}
\proof 
By Proposition \ref{rich}, the endofunctor $\GG$ is cartesian, and the monad $(\GG,\sigma,\pi)$ is strongly $\P$-cartesian with a $\P$-cartesian $\sigma$; so, $\sigma$ is the equalizer of the pair $(\sigma_\GG,\GG(\sigma))$, and the monad is half-cartesian. Then apply Corollary \ref{auton}: accordingly the adjoint pair $(U^\GG,F^\GG): Grd\EE \splito \EE$ is autonomous. As for the second point, starting with an internal groupoid $X_{\bullet}$, observe that the endofunctor $Dec$ on $Cat\EE$ is stable on $Grd\EE$:
$$\xymatrix@=12pt{
	DecX_{\bullet}\ar[dd]_{\epsilon_{X_{\bullet}}}: & R_3[d_0]\ar[dd]_{d_4} \ar@<-3ex>[rr]_{d_3} \ar@<-1ex>[rr]|{d_2} \ar@<1ex>[rr]|{d_1} \ar@<3ex>[rr]^{d_0}	&& R_2[d_0] \ar[dd]_{d_3}\ar@<-2ex>[rr]_{d_2} \ar[rr]|{d_1} \ar@<2ex>[rr]^{d_0}  &&	R[d_0] \ar[dd]^{d_2}\ar@<-2ex>[rr]_{d_1} 
	\ar@<2ex>[rr]^{d_0}  && X_1\ar[ll]|{s_0} \ar[dd]^{d_1}\\
	&&&&\\
	X_{\bullet}: & R_2[d_0] \ar@<-3ex>[rr]_{d_3} \ar@<-1ex>[rr]|{d_2} \ar@<1ex>[rr]|{d_1} \ar@<3ex>[rr]^{d_0}	&&	R[d_0]\ar@<-2ex>[rr]_{d_2} \ar[rr]|{d_1} \ar@<2ex>[rr]^{d_0}  &&	X_1 \ar@<-2ex>[rr]_{d_1} 
	\ar@<2ex>[rr]^{d_0}  && X_0\ar[ll]|{s_0}
}
$$
since the upper horizontal diagram is the groupoid $R_{\bullet}[d_0]$ and that this $Dec X_{\bullet}$ is nothing but $F^\GG.U^\GG(X_{\bullet})$.

As for the last point, consider the following vertical discrete fibration in $Cat\EE$:
$$\xymatrix@=15pt{
	DecX_{\bullet}\ar[dd]_{\epsilon_{X_{\bullet}}}: & X_4\ar[dd]_{d_4} \ar@<-3ex>[rr]_{d_3} \ar@<-1ex>[rr]|{d_2} \ar@<1ex>[rr]|{d_1} \ar@<3ex>[rr]^{d_0}	&& X_3 \ar[dd]_{d_3}\ar@<-2ex>[rr]_{d_2} \ar[rr]|{d_1} \ar@<2ex>[rr]^{d_0}  &&	X_2 \ar[dd]^{d_2}\ar@<-2ex>[rr]_{d_1} 
	\ar@<2ex>[rr]^{d_0}  && X_1\ar[ll]|{s_0} \ar[dd]^{d_1}\\
	&&&&\\
	X_{\bullet}: & X_3 \ar@<-3ex>[rr]_{d_3} \ar@<-1ex>[rr]|{d_2} \ar@<1ex>[rr]|{d_1} \ar@<3ex>[rr]^{d_0}	&&	X_2\ar@<-2ex>[rr]_{d_2} \ar[rr]|{d_1} \ar@<2ex>[rr]^{d_0}  &&	X_1 \ar@<-2ex>[rr]_{d_1} 
	\ar@<2ex>[rr]^{d_0}  && X_0\ar[ll]|{s_0}
}
$$
and suppose the upper horizontal part is groupoid. Since the vertical discrete fibration $h_{\bullet}$ has its $h_0=d_1$ split by $s_0$, then the lower horizontal row is a groupoid: denote by $\gamma: X_2\to X_2$ the involutive mapping producing the inversion in the upper groupoid. Then $d_2.\gamma.s_1:X_1\to X_1$ determines an involutive map which produces an inversion for the lower row and makes it a groupoid as well. 
\endproof
\begin{coro}\label{KlG}
	The category $\Pt\EE$ is a subcategory of $Kl\GG$; the morphisms of this last category are the commutative squares in $\EE$:
	$$\xymatrix@=20pt{
		X \ar@<-4pt>[d]_{g} \ar[r]^-{x} & X'  \ar@<-4pt>[d]_{g'}  \\
		Y \ar[r]_{y} \ar@{ >.>}[u]_t & Y' \ar@{ >.>}[u]_{t'} }
	$$
	not necessarily respecting the sections.
\end{coro}
\proof
The first assertion is a consequence of Proposition \ref{incl}. The morphisms in $Kl\GG$ between the two vertical objects of the previous diagram are given by the internal functors between $R[g]$ and $R[g']$ which produce (and are produced by) the commutative diagrams in question since $g$ and $g'$, being split, are the quotients of these kernel equivalence relations.
\endproof

\subsection{The monad $(T_{X_{\bullet}},\lambda_{X_{\bullet}},\mu_{X_{\bullet}})$ when $X_{\bullet}$ is a groupoid}\label{truc}  
\begin{prop}
Given any internal category $X_{\bullet}$, the monad $(T_{X_{\bullet}},\lambda_{X_{\bullet}},\mu_{X_{\bullet}})$ is hypercartesian if and only if $X_{\bullet}$ is a groupoid. 	
\end{prop}
\proof
Suppose $X_{\bullet}$ is a groupoid. We have to show that the monad is hypercartesian. For that, let us reproduce below a part of the upper part of the diagram of Section \ref{T_X} and let us complete it on the left hand side:
$$\xymatrix@=25pt{
	d_1^*(Z)  \ar@<-2pt>[d]^{d_1^*(h)} & (d_1.d_2)^*(Z)\ar@<1ex>[l]^{\delta_2^h} \ar@<-1ex>@{.>}[l]_{\delta_1^h}\ar@<-3pt>[d]^{(d_1.d_2)^*(h)} && (d_1.d_2.d_3)^*(Z) \ar@<-3pt>[d]^{(d_1.d_2.d_3)^*(h)}\ar@<2ex>[ll]^{\delta_3^h} \ar[ll]|{\delta_2^h} \ar@<-2ex>@{.>}[ll]_{\delta_1^h}  \\
	X_1  & R[d_1]\ar@<1ex>[l]^{d_2}\ar@<-1ex>@{.>}[l]_{d_1} && \ar@<2ex>[ll]^{d_3} \ar[ll]|{d_2} \ar@<-2ex>@{.>}[ll]_{d_1} R_2[d_1]
}
$$
In this diagram, any commutative vertical square is a pullback. Since $X_{\bullet}$ is a groupoid, we get the following kernel equivalence relations on the lower row:
$$\xymatrix@=25pt{
	X_1   & R[d_1] \ar[l]_{d_1} && R_2[d_1]  \ar@<1ex>[ll]^{d_2} \ar@<-1ex>[ll]_{d_1}  
}
$$
which is lifted by pullback on the upper row as a kernel equivalence relation:
$$\xymatrix@=25pt{
	d_1^*(Z)   & (d_1.d_2)^*(Z) \ar[l]_{\delta_1^h} && (d_1.d_2.d_3)^*(Z)  \ar@<1ex>[ll]^{\delta_2^h} \ar@<-1ex>[ll]_{\delta_1^h}  
}
$$
which is the hypercartesian condition for the cartesian monad $(T_{X_{\bullet}},\lambda_{X_{\bullet}},\mu_{X_{\bullet}})$.

Conversely, suppose this monad is hypercartesian. The hypercartesian condition applied to the terminal object $1_{X_0}$ of $\EE/X_0$ says that in the following diagram which is nothing but the $DecX^{op}_{\bullet}$:
$$\xymatrix@=25pt{
	X_1 \ar[r]|{s_1} & X_2\ar@<1ex>[l]^{d_2}\ar@<-1ex>@{.>}[l]_{d_1} && X_3\ar@<2ex>[ll]^{d_3} \ar@{.>}[ll]|{d_2} \ar@<-2ex>@{.>}[ll]_{d_1} 
}
$$
the dotted part of the diagram is a kernel equivalence relation, and consequently that  $DecX^{op}_{\bullet}$ is a groupoid. According to Corollary \ref{GGaut}, $X^{op}_{\bullet}$, and thus $X_{\bullet}$, is a groupoid. 
\endproof
This gives a conceptual way to a straightforward result in $Set$:
\begin{coro}\label{fib=cofib}
 The domain of any discrete fibration above a  groupoid in $Cat\EE$ is a groupoid as well.
\end{coro}
\proof
We know that any $T_{X_{\bullet}}$-algebra produces a discrete fibration above the groupoid $X_{\bullet}$:
$$\xymatrix@=10pt{
	 Z_2 \ar[dd]_{h_2}\ar@<-2ex>[rr]_{d_2} \ar[rr]|{d_1} \ar@<2ex>[rr]^{d_0}  &&	Z_1 \ar[dd]^{h_1}\ar@<-2ex>[rr]_{d_1} 
	\ar@<2ex>[rr]^{d_0}  && Z_0\ar[ll]|{s_0} \ar[dd]^{h_0}\\
	&&&&\\
	R[d_0]\ar@<-2ex>[rr]_{d_2} \ar[rr]|{d_1} \ar@<2ex>[rr]^{d_0}  &&	X_1 \ar@<-2ex>[rr]_{d_1} 
	\ar@<2ex>[rr]^{d_0}  && X_0\ar[ll]|{s_0}
}
$$
By Proposition \ref{cartadj}, any $T_{X_{\bullet}}$-algebra $d_0: T_{X_{\bullet}}(h)\to h$ produces a kernel equivalence relation:
$$
\xymatrix@C=45pt{
	T_{X_{\bullet}}^2(h)  \ar@<-7pt>[r]_-{\mu_{X_{\bullet}h}}
	\ar@<7pt>[r]^-{T_{X_{\bullet}}(d_0)} &
	T_{X_{\bullet}}(h) \ar@{->>}[r]^-{d_0} \ar[l]|-{} & h }
$$
which, here, is nothing but:
$$
\xymatrix@C=45pt{
	Z_2  \ar@<7pt>[r]^-{d_1}
	\ar@<-7pt>[r]_-{d_0} &
	Z_1 \ar@{->>}[r]^-{d_0} \ar[l]|-{} & Z_0 }
$$
and shows that the category  $Z_{\bullet}$ is a groupoid.
\endproof 

\section{$T$-categories}\label{Tcat}

Now, let $(T,\lambda,\mu)$ be any monad on $\EE$. The notion of $T$-category has been introduced by A. Burroni in \cite{AB} as a mix of a relational algebra in the sense of Barr \cite{Barr1} and of something which looks like a kind of internal category, but shifted by this monad. He first introduced the notion of \emph{pointed $T$-graph} in $\EE$ as a triple $(d_0,\delta_1,s_0)$ of maps:
	$$\xymatrix@=16pt{
	& X_1 \ar[dl]_{d_0} \ar[dr]^{\delta_1}  \\
	X_0  \ar@<-1ex>[ur]_{s_0} \ar[rr]_{\lambda_{X_0}}  && T(X_0)
}
$$
such that $d_0.s_0=1_{X_0}$ and $\delta_1.s_0=\lambda_{X_0}$ (Axioms 1). According to our notations, it is nothing but a reflexive graph in the Kleisli category $KlT$:
	$$\xymatrix@=16pt{
	X_1 \ar@<-2ex>[rr]_{\bar F^T(d_0)} 
  \ar@<2ex>@{.>}[rr]^{"\delta_1"}   && X_0 \ar[ll]|{\bar F^T(s_0)}
}
$$
Then building the pullback of $\delta_1$ along $T(d_0)$ in $\EE$ (in plain arrows in the following diagram):
	$$\xymatrix@=16pt{
		&& X_2  \ar@<-1ex>[dl]_{d_0^1} \ar@{-->}@<1ex>[dl]^{d_1^1} \ar[dr]^{\delta_2} && \\
	& X_1 \ar[dl]_{d_0} \ar[dr]^{\delta_1} \ar[ur]|{s_0^1} && T(X_1) \ar@<-1ex>[dl]_{T(d_0)}  \ar[dr]^{T(\delta_1)}\\
	X_0  \ar@<-1ex>[ur]_{s_0} \ar@{.>}[rr]_{\lambda_{X_0}}  && T(X_0) \ar[ur]_{T(s_0)} && T^2(X_0) \ar[ll]^{\mu_{X_0}} 
}
$$
Burroni first observed that the section $s_0$ of $d_0$ produces, through the existence of $T(s_0)$ a section $s_0^1$ of $d_0^1$ such that $\delta_2.s_0^1=T(s_0).\delta_1$ (Observation 2) and that the identity $\delta_1.s_0.d_0=\lambda{X_0}.d_0=T(d_0).\lambda_{X_1}$ produces a map $s_1^1:X_1\to X_2$ such that $d_0^1.s_1^1=s_0.d_0$ and $\delta_2.s_1^1=\lambda X_1$  (Observation 3).

Then \emph{he demanded a "composition" map} $d_1^1: X_2\to X_1$ in $\EE$ satisfying: $d_0.d_1^1=d_0.d_0^1$ and $\delta_1.d_1^1=\mu_{X_0}.T(\delta_1).\delta_2$ (Axioms 4) which, with our notation, delineates the beginning of a $2$-truncated simplicial object in $KlT$:
	$$\xymatrix@=18pt{
	X_2\ar@<-4ex>[rr]_{\bar F^T(d_0^1)} \ar[rr]|{\bar F^T(d_1^1)} \ar@<4ex>@{.>}[rr]^{"\delta_2"}   &&	X_1 \ar@<-2ex>[rr]_{\bar F^T(d_0)} 
	\ar@<2ex>@{.>}[rr]^{"\delta_1"} \ar@<2ex>[ll]|{\bar F^T(s_0^1)} \ar@<-2ex>[ll]|{\bar F^T(s_1^1)} && X_0 \ar[ll]|{\bar F^T(s_0)}
}
$$
Then Buronni demands Axioms 7 (\emph{neutrality}): $d_1^1.s_0^1=1_{X_1}$ and $d_1^1.s_1^1=1_{X_1}$ which completes the previous diagram into a plain $2$-truncated simplicial object in $KlT$.
Finally, constructing the pullback of $\delta_2$ along $T(d_0^1)$:
$$\xymatrix@=16pt{
	&& X_3  \ar@<-3ex>[dl]_{d_0^2} \ar@{-->}[dl]|<<<<{d_1^2} \ar@<2ex>@{-->}[dl]^{d_2^2} \ar[dr]^{\delta_3} && \\
	& X_2 \ar@<1ex>[dl]^{d_1^1}  \ar@<-1ex>[dl]_{d_0^1} \ar[dr]^{\delta_2} \ar@<2ex>[ur]|{s_0^2} && T(X_2) \ar@<-2ex>[dl]_{T(d_0^1)} \ar@<2ex>[dl]^<<<<{T(d_1^1)} \ar[dr]^{T(\delta_2)}\\
	X_1 \ar[ur]|{s_0^1}   && T(X_1)\ar[ur]|{T(s_0^1)} && T^2(X_1) \ar[ll]^{\mu_{X_1}}  
}
$$
Burroni observed that:\\
1) from: $\delta_1.d_0^1.d_0^2=T(d_0).\delta_2.d_0^2=T(d_0).T(d_0^1).\delta_3=T(d_0).T(d_1^1).\delta_3$\\ we get a morphism $d_1^2: X_3\to X_2$ such that $d_0^1.d_1^2=d_0^1.d_0^2$ and $\delta_2.d_1^2=T(d_1^1).\delta_3$ (Observations 5);\\ 2) and from:
$\delta_1.d_1^1.d_0^2=\mu_{X_0}.T(\delta_1).\delta_2.d_0^2=\mu_{X_0}.T(\delta_1).T(d_0^1).\delta_3$\\$=\mu_{X_0}.T^2(d_0).T(\delta_2).\delta_3=T(d_0).\mu_{X_1}.T(\delta_2).\delta_3$, we get a morphism $d_2^2: X_3\to X_2$ such that: $d_0^1.d_2^2=d_1^1.d_0^2$ and $\delta_2.d_2^2=\mu_{X_1}.T(\delta_2).\delta_3$ (Observations 6).\\
Then he added Axiom 8 (\emph{associativity}) $d_1^1.d_1^2=d_1^1.d_2^2$ to complete the definition of a $T$-category.

Further observations: the splitting $s_0^1$ of $d_0^1$ produces, via the map $T(s_0^1)$, a splitting $s_0^2: X_1\to X_2$ of $d_0^2$ such that $\delta_3.s_0^2=T(s_0^1).\delta_2$. And from: $d_1^1.s_0^1.d_0=\lambda_{X_0}.d_0=T(d_0).\lambda_{X_1}$ we get a map $s_1^1: X_1\to X_2$ such that $d_0^1.s_1^1=s_0.d_0$ and $\delta_2.s_1^1=\lambda_{X_1}$. Finally from $\delta_2.s_1^1.d_0^1=\lambda_{X_1}.d_0^1=T(d_0^1).\lambda_{X_1}$, we get a map $s_2^2:X_2\to X_3$ such that $d_0^2.s_2^2=s_1^1.d_0^1$ and $\delta_3.s_2^2=\lambda_{X_2}$. Again, with our notation, this delineates a  $3$-truncated object in the category $KlT$ (where the higher degeneracies are omitted):
$$\xymatrix@=18pt{
	X_3 \ar@<-7ex>[rr]_{\bar F^T(d_0^2)} \ar@<7ex>@{.>}[rr]^{"\delta_3"} \ar@<-3ex>[rr]|{\bar F^T(d_1^2)} \ar@<3ex>[rr]|{\bar F^T(d_2^2)}&& X_2\ar@<-4ex>[rr]_{\bar F^T(d_0^1)} \ar[rr]|{\bar F^T(d_1^1)} \ar@<4ex>@{.>}[rr]^{"\delta_2"} \ar[ll]|{\bar F^T(s_1^2)} \ar@<5ex>[ll]|{\bar F^T(s_0^2)} \ar@<-5ex>[ll]|{\bar F^T(s_2^2)} &&	X_1 \ar@<-2ex>[rr]_{\bar F^T(d_0)} 
	\ar@<2ex>@{.>}[rr]^{"\delta_1"} \ar@<2ex>[ll]|{\bar F^T(s_0^1)} \ar@<-2ex>[ll]|{\bar F^T(s_1^1)} && X_0 \ar[ll]|{\bar F^T(s_0)}
}
$$

A morphism of $T$-categories, namely a \emph{$T$-functor}, is a morphism $(f_0,f_1)$ of pointed $T$-graph:
$$\xymatrix@=7pt{
	&& T(X_0) \ar[rrr]^{T(f_0)} &&& T(Y_1)\\
	X_1 \ar[drr]_{d_0} \ar[urr]^{\delta_1}\ar[rrr]^{f_1} &&&	 Y_1\ar[drr]^{d_0} \ar[urr]_{\delta_1}\\
	&& X_0 \ar[rrr]_{f_0} &&&	 Y_0 
}
$$
which preserves the "composition maps" $d_1^1$. Any  $T$-functor naturally induces a morphism of $3$-truncated simplicial objects in $KlT$.
Whence the category $T$-Cat$\EE$ of $T$-categories whose objects will be denoted  $X^T_{\bullet}$ and morphisms $h^T_{\bullet}: X^T_{\bullet}\to Y^T_{\bullet}$. There is a forgetful functor $(\;)_T$:$T$-Cat$\EE \to \EE$ associating the object $X_0$ with the whole $T$-category $X_0 \stackrel{d_0}{\leftarrow} X_1 \stackrel{\delta_1}{\rightarrow} T(X_0)$.

We get an injective fully faithful injective functor $\TC: \AlgT \into {\rm T}$-Cat$\EE$ where $\TC_{\EE}(X,\xi)$ has the following underlying pointed $T$-graph:
$$\xymatrix@=12pt{
	& T(X) \ar[dl]_{\xi} \ar[dr]^{1_{T(X)}} &\\
	X \ar[rr]_{\lambda_X} \ar@<-1ex>[ru]_>>{\lambda_X} && T(X)
}
$$
the structure of $T$-category being produced by the following diagram:
$$\xymatrix@=12pt{
	&& T^2(X)  \ar@<-1ex>[dl]_{T(\xi)} \ar@{-->}[dl]^{\mu_X} \ar[dr]^{1_{T^2(X)}} && \\
	& T(X) \ar[dl]_{\xi} \ar[dr]^{1_{T(X)}} && T^2(X) \ar[dl]_{T(\xi)}  \ar[dr]^{T(1_{T(X)})}\\
	X  \ar@<-1ex>[ur]_{\lambda_X} \ar@{.>}[rr]_{\lambda_{X}}  && T(X) && T^2(X) \ar@{.>}[ll]^{\mu_{X}} 
}
$$
The functor $\TC_{\EE}$ makes commute the following diagram:
$$\xymatrix@=25pt{
	\AlgT  \ar@{ >->}[r]^{\TC_{\EE}} \ar@<-2pt>[d]_{U^T} & {\rm T-}Cat\EE  \ar@<2pt>[d]^{(\;)_T} \\
	\EE   \ar@{=}[r]  & \EE \
}
$$
In this way, the category ${\rm T}$-Cat$\EE$ appears as a kind of natural extension of the category $\AlgT$. 

\noindent\textbf{Warning:} However an internal category $\xymatrix@=8pt{
X_{\bullet}:	X_1 \ar@<-1ex>[rr]_{d_0} 
	\ar@<1ex>[rr]^{d_1}   && X_0 \ar[ll]|<<<<{s_0}
}
$ in $\EE$ does not induce in general a structure of $T$-category on the following pointed $T$-graph:
$$\xymatrix@=12pt{
	& X_1 \ar[dl]_{d_0} \ar[dr]^{\lambda_{X_0}.d_1}  \\
	X_0  \ar@<-1ex>[ur]_{s_0} \ar[rr]_{\lambda_{X_0}}  && T(X_0)
}
$$
since the pullback of $T(d_0)$ along $\lambda_{X_0}.d_1$ does not coincide with $X_2$ in general. When the endofunctor $T$ is cartesian, so is $\bar F^T: \EE\to KlT$. However, for the same reason, the internal category $\bar F^T(X_{\bullet})$ in $KlT$ does not coincide with a $T$-category.

\begin{prop}
	Suppose $T$ and $\lambda$ cartesian. Then the image by the inclusion $\bar F^T:\EE\into KlT$ of an internal category is a $T$-category. So, we really get an inclusion functor $Cat(\bar F^T): Cat\EE \into T$-$Cat\EE$.
\end{prop}
\proof
Start with an internal category $X_{\bullet}$ in $\EE$. We have $\lambda_{X_0}.d_1=T(d_1).\lambda_{X_1}$. Now, since the following whole rectangle is a pullback in $\EE$
$$\xymatrix@=5pt{
	X_2\ar[rr]^{\lambda_{X_2}} \ar[dd]_{d_0} &&	T(X_2) \ar[dd]_{T(d_0)} \ar[rr]^{T(d_2)}  && T(X_1) \ar[dd]^{T(d_0)} && \\
	\\
	X_1 \ar[rr]_{\lambda_{X_1}} &&	T(X_1) \ar[rr]_{T(d_1)}  && T(X_0) 
}
$$
the composition map $d_1^1:X_2\to X_1$ of $X_{\bullet}$  in $\EE$ produces the map $d_1^1$ demanded by the definition of a $T$-category. The satisfaction of the other axioms immedialely follows.
\endproof

\section{When $T$-categories in $\EE$ coincide with (a special kind of) internal categories in $KlT$}

 In this section, we are going to investigate the $T$-categories in the  setting of the $\Sigma$-cartesian monads and to show that some specific class of $T$-categories coincides with some specific class of internal categories in the Kleisli category $KlT$ of the monad. This will be the case thanks to the following:
\begin{prop}
	Given any pullback stable class $\Sigma$ in $\EE$ and any $\Sigma$-cartesian monad $(T,\lambda,\mu)$, then:\\
	1) the class $\Sigma$ is pullback stable in $KlT$;\\
	2) given a pair $(g,h)$ of maps in $\Sigma\times\EE$ with same codomain $Z$, if there is a map $\phi:X\cdots>Y$ in $KlT$ such that $g.\phi=h$, then $\phi$ belongs to $\EE$.
	
\end{prop}
\proof
First, since $\lambda$ is the equalizer of the pair $(\lambda_T,T(\lambda))$, $\EE$ is a subcategory of $KlT$.
1) By the same proof as the one of Proposition \ref{spullback}, restricted to the maps in $\Sigma$, any map in $\Sigma$ is pullback stable in $KlT$.\\
2) Again, by a careful inspection, the same proof as the one of Proposition \ref{c3from2} works, when it is restricted to the maps in $\Sigma$.
\endproof 

\subsection{$T$-categories and $\Sigma$-cartesian monads}
 
Now, let be given  any pointed $T$-graph in $\EE$:
$$\xymatrix@=16pt{
	& X_1 \ar[dl]_{d_0} \ar[dr]^{\delta_1}  \\
	X_0  \ar@<-1ex>[ur]_{s_0} \ar[rr]_{\lambda_{X_0}}  && T(X_0)
}
$$
\begin{prop}\label{ssigma}
Let $\Sigma$ be a pullback stable class of morphisms in $\EE$ and $(T,\lambda,\mu)$ be a $\Sigma$-cartesian monad. Then there is a bijection between the $T$-categories  $X^T_{\bullet}$ having its underlying pointed $T$-graph with leg $d_0 \in \Sigma$ and the internal categories in $KlT$:
$$\xymatrix@=18pt{
 X_{\bullet}:&&	X_3 \ar@<-7ex>[rr]_{d_0^2} \ar@<7ex>@{.>}[rr]^{"\delta_3"} \ar@<-3ex>[rr]|{d_1^2} \ar@<3ex>[rr]|{d_2^2}&& X_2\ar@<-4ex>[rr]_{d_0^1} \ar[rr]|{d_1^1} \ar@<4ex>@{.>}[rr]^{"\delta_2"} \ar[ll]|{s_1^2} \ar@<5ex>[ll]|{s_0^2} \ar@<-5ex>[ll]|{s_2^2} &&	X_1 \ar@<-2ex>[rr]_{d_0} 
	\ar@<2ex>@{.>}[rr]^{"\delta_1"} \ar@<2ex>[ll]|{s_0^1} \ar@<-2ex>[ll]|{s_1^1} && X_0 \ar[ll]|{s_0}
}
$$
$KlT$ having leg $d_0: X_1\to X_0 \in \Sigma$ and section $s_0\in \EE$.
\end{prop}
\proof Starting with such a $T$-category $X^T_{\bullet}$, according to the previous proposition, the following diagrams are pullbacks  in $KlT$  since $d_0$ is in $\Sigma$:
$$\xymatrix@=7pt{
	X_3\ar[rr]^{d_0^2} \ar@{.>}[dd]_{"\delta_3"} &&	X_2 \ar@{.>}[dd]_{"\delta_2"} \ar[rr]^{d_0^1}  && X_1 \ar@{.>}[dd]^{"\delta_1"} && \\
	\\
	X_2 \ar[rr]_{d_0^1} &&	X_1 \ar[rr]_{d_0}  && X_0 
}
$$
and the above $3$-truncated simplicial object in $KlT$ is underlying an internal category in $KlT$.

Conversely, starting with any internal category in $KlT$ with $d_0:X_1\to X_0 \in \Sigma$ and $s_0\in \EE$, the pullbacks involved in the definition of an internal category are obtained as above, and then they coincide with Burroni's construction. By the second part of the same proposition and the identity $d_0.d_1^1=d_0.d_0^1$, the composition map $d_1^1:X_2\to X_1$ necessarily belongs to $\EE$, since $d_0\in \Sigma$; in the same way, since $d_0^1\in \Sigma$ as well, by $d_0^1.d_1^2=d_0^1.d_0^2$, the map $d_1^2$ is in $\EE$ too; finally the identity $d_0^1.d_2^2=d_1^1.d_0^2$ assures us that $d_2^2$ is in $\EE$; accordingly, the internal category $X_{\bullet}$ in $KlT$ is underlying a $T$-category in $\EE$.
\endproof 
Similarly, by the same previous proposition, we know that any internal functor $(h_0,h_1)$  between the internal categories in $KlT$ with $d_0\in \Sigma$ having $h_0\in \EE$ is such that $h_1\in \EE$; accordingly the full subcategory  $T_{\Sigma}$-$Cat\EE$ of $T$-$Cat\EE$ whose objects are the $T$-categories with $d_0\in \Sigma$ coincides with the  subcategory $Cat_{\Sigma}KlT$ of internal categories in $KlT$ having $d_0\in\Sigma$ and $s_0\in\EE$ and internal functors $(h_0,h_1)$ in $KlT$ with $h_0\in \EE$. It is obtained by the following pullback:  
$$\xymatrix@=20pt{
	{\rm T_{\Sigma}-}Cat\EE  \ar[rr] \ar@<-2pt>[d]_{} && Cat\AlgT  \ar@<2pt>[d]^{D^0_{\AlgT}}\\
	Pt_{\Sigma}\EE   \ar@{ >->}[r]  & Pt\EE   \ar[r]_{Pt(F^T)}  & Pt\AlgT
}
$$
where $Pt_{\Sigma}\EE$ denotes the full subcategory of $Pt\EE$ whose objets are the split epimorphism $(f,s)$ in $\EE$ with $f\in \Sigma$.

When, in addition, \emph{the class $\Sigma$ contains the identity maps, is stable under composition and left cancellable}, the situation becomes even clearer since we are now assured that the map $s_0$ and the composition map $d_1:X_2\to X_1$ belong to $\Sigma$. So, we get a fully faithful inclusion $j: Cat_{\Sigma}\EE\into$ $T_{\Sigma}$-$Cat\EE$ where $Cat_{\Sigma}\EE$ denotes the  full subcategory of $Cat\EE$ whose objects are the internal categories in $\Sigma$.

\subsection{$T$-categories and cartesian monads}

In this way, we get our more meaningful result:
\begin{theo}\label{mmain}
	When the monad $(T,\lambda,\mu)$ is cartesian, a $T$-category coincides with an internal category $X^T_{\bullet}$ in $KlT$ whose leg  $d_0$ belongs to the subcategory $\EE$. A $T$-functor, coincides with an internal functor in $KlT$ whose image by the functor $(\;)_0: Cat(KlT)\to KlT$ belongs to the subcategory $\EE$.  The image by the cartesian inclusion functor $\bar F^T: \EE\into KlT$ of any internal category is a $T$-category. The image of any $T$-category by the fully faithful functor $K_T:KlT\to \AlgT$ produces an internal category in $\AlgT$.
	\end{theo}
\proof
It is a corollary of the previous proposition where $\Sigma=\EE$.
\endproof
By Proposition III.2.21 in \cite{AB}, Burroni observed that, when the monad $(T,\lambda,\mu)$ is cartesian, the image by the functor $\bar U^T=U^T.K_T: KlT\to \EE$ of the $3$-truncated simplicial object in $KlT$ induced by a $T$-category in $\EE$ is an internal category in $\EE$, but he did not produced the previous characterization; for that the Proposition \ref{c3from2} concerning the behaviour of the maps  of $\EE$ inside $KlT$ and the Proposition \ref{spullback} concerning the existence of a certain class of pullbacks in $KlT$ are needed.
According to the previous proposition, the category $T$-$Cat\EE$ is defined by any of the following pullbacks:
$$\xymatrix@=15pt{
	{\rm T-}Cat\EE  \ar@{ >->}[r] \ar@<-2pt>[d]_{} & CatKlT  \ar@<2pt>[d]^{D^0_{KlT}} && 	{\rm T-}Cat\EE  \ar[r] \ar@<-2pt>[d]_{} & Cat\AlgT  \ar@<2pt>[d]^{D^0_{\AlgT}}\\
	\Pt\EE   \ar@{ >->}[r]_{\Pt(\bar F^T)}  & \Pt KlT && \Pt\EE   \ar[r]_{\Pt(F^T)}  & \Pt\AlgT
}
$$
 Let us denote by $(\;)_T: T$-$Cat\EE \to \EE$ the forgetful functor associating with any $T$-category $X^T_{\bullet}$ its "object of objects" $X_0$.
\begin{prop}
	When the monad $(T,\lambda,\mu)$ is cartesian, the category $T$-$Cat\EE$ has pullbacks and
	the forgetful functor $(\;)_T$ is cartesian; it is a fibration such that, in the following commutative diagram, the inclusion $Cat(\bar F^T):Cat\EE \into T$-$Cat\EE$ is fully faithful, cartesian and preserves the cartesian maps:
	$$\xymatrix@=25pt{
		{Cat\EE\;}  \ar@{ >->}[r]^{Cat(\bar F^T)} \ar@<-2pt>[d]_{(\;)_0} & {\rm T-}Cat\EE  \ar@<2pt>[d]^{(\;)_T} \\
		\EE   \ar@{=}[r]  & \EE \
	}
	$$
\end{prop}
\proof
Let $f_{\bullet}: X^T_{\bullet}\to Z^T_{\bullet}$ and $g_{\bullet}: Y^T_{\bullet}\to Z^T_{\bullet}$ be a pair of $T$-functors. Consider the following levelwise pullbacks in $\EE$:
$$\xymatrix@=7pt{
	P_0\ar[rr]^{p_0^X} \ar[dd]_{p_0^Y} &&	X_0 \ar[dd]^{f_0} && P_1 \ar[dd]_{p_1^Y} \ar[rr]^{p_1^X}  && X_1 \ar[dd]^{f_1} && \\
	\\
	Y_0 \ar[rr]_{g_0} && Z_0 && Y_1 \ar[rr]_{g_1}  && Z_1 
}
$$
The split epimorphisms $(d_0,s_0)$ in $\EE$ produce a  split epimorphism $(d_0^P,s_0^P): P_1\splito P_0$ in $\EE$. And since the injection $\EE\into KlT$ is cartesian (Proposition \ref{onT}), the maps $"\delta_1": X_1\cdots> X_0$ in $KlT$ produces a map $"\delta_1":P_1\cdots> P_0$ in $KlT$; from that the structure $P^T_{\bullet}$ of $T$-category on this induced pointed graph follows. By this construction,  the functors $(\;)_T$ and $Cat(\bar F^T)$ are cartesian.

To show that this functor is a fibration, we have first to check that the classical construction of the cartesian maps above a map $f: X\to Y_0$ is valid in $KlT$, namely to build some joint pullbacks in $KlT$. So let $Y^T_{\bullet}$ be a $T$-category and $f:X\to Y_0$ any map in $\EE$. The following diagram where any square is a pullback in $KlT$ makes it explicit:
$$\xymatrix@=4pt{
	X_1 \ar[dd]_{\bar{\phi}} \ar[rr]^{\check{\bar f}}	&& \check X \ar@{.>}[rr]^{"\check\delta_1"} \ar[dd]_{\phi} &&	X \ar[dd]^{f}   \\
	\\
	\bar X	\ar[dd]_{\bar d_0}\ar[rr]_{\bar f} && Y_1 \ar@{.>}[rr]_{"\delta_1"} \ar[dd]^{d_0} && Y_0    \\
	\\
	X \ar[rr]_{f} && Y_0
}
$$
Whence a morphism of pointed graphs in $KlT$: 
$$\xymatrix@=7pt{
	X_1\ar@<-4pt>[dd]_{\bar d_0.\bar{\phi}} \ar@<4pt>@{.>}[dd]^{"\check\delta _1".\check{\bar f}} \ar[rrrrr]^{\bar f.\bar{\phi}}  &&&&& Y_1  \ar@<-4pt>[dd]_{d_0} \ar@<4pt>@{.>}[dd]^{"\delta_1"} \\
	&&&&&&&\\
	X \ar@{ >->}[uu]|{} \ar[rrrrr]_{f}  
	&&&&& Y_0     \ar@{ >->}[uu]
}
$$
which, by general arguments, endows  the left hand side reflexive graph with an internal category structure $X^T_{\bullet}$ in $KlT$. By Theorem \ref{mmain}, this internal category in $KlT$ is a $T$-category.

We have now to check the universal property: so, let $Z^T_{\bullet}$ be another $T$-category and $g_{\bullet}: Z^T_{\bullet}\to Y^T_{\bullet}$ a $T$-functor, such that $g_0=f.h$ for some $h:Z_0\to X$ in $\EE$. Again, by general arguments, we certainly get an internal functor $h_{\bullet}: Z^T_{\bullet}\to X^T_{\bullet}$ in $KlT$ such that $h_0=h\in \EE$. Since both $Z_{\bullet}$ and $X_{\bullet}$ are $T$-categories and $h_0$ belongs to $\EE$, then $h_{\bullet}$ is a $T$-functor, see Theorem \ref{mmain}. Since $\bar F^T$ preserves pullbacks, $Cat(\bar F^T)$ preserves the fully faithful internal functors, namely the cartesian maps with respect to $(\;)_0$.
\endproof
In our context, the construction of the endofunctor $Dec$ can be extended to $T$-categories:
\begin{prop}
	Given any cartesian monad $(T,\lambda,\mu)$, there is an endofunctor $Dec$ on $T$-$Cat\EE$ which mimicks the endofunctor $Dec$ on $Cat\EE$. However the co-unit $\epsilon$ does not belong to $T$-$Cat\EE$.	
\end{prop}
\proof
Consider the upper part of the following vertical diagram in $KlT$:
$$\xymatrix@=15pt{
	DecX_{\bullet}\ar[dd]_{\epsilon_{X_{\bullet}}}: & X_4\ar@{.>}[dd]_{\delta_4} \ar@<-3ex>[rr]_{d_3} \ar@<-1ex>[rr]|{d_2} \ar@<1ex>[rr]|{d_1} \ar@<3ex>[rr]^{d_0}	&& X_3 \ar@{.>}[dd]_{"\delta_3"}\ar@<-2ex>[rr]_{d_2} \ar[rr]|{d_1} \ar@<2ex>[rr]^{d_0}  &&	X_2 \ar@{.>}[dd]^{"\delta_2"}\ar@<-2ex>[rr]_{d_1} 
	\ar@<2ex>[rr]^{d_0}  && X_1\ar[ll]|{s_0} \ar@{.>}[dd]^{"\delta_1"}\\
	&&&&\\
	X_{\bullet}: & X_3 \ar@<-3ex>@{.>}[rr]_{"\delta_3"} \ar@<-1ex>[rr]|{d_2} \ar@<1ex>[rr]|{d_1} \ar@<3ex>[rr]^{d_0}	&&	X_2\ar@<-2ex>@{.>}[rr]_{"\delta_2"} \ar[rr]|{d_1} \ar@<2ex>[rr]^{d_0}  &&	X_1 \ar@<-2ex>@{.>}[rr]_{"\delta_1"} 
	\ar@<2ex>[rr]^{d_0}  && X_0\ar[ll]|{s_0}
}
$$
where $X_4$ is defined by the pullback in $KlT$ of the map $d_0:X_3\to X_2\in \EE$ along the map $\delta_3:X_3\cdots>X_2\in KlT$.
\endproof

\begin{prop}\label{Rad}
	Given any cartesian monad $(T,\lambda,\mu)$, the fully faithful inclusion $Cat(\bar F^T): Cat\EE \into  T$-$Cat\EE$ admit a right adjoint $\RR$ which preserves the cartesian maps (=fully faithful functors).
\end{prop}
\proof
Let us start with a $T$-category $X^T_{\bullet}$ and define $\RR(X^T_{\bullet})$ by the following fully faithful internal functor in $Cat\EE$:
$$\xymatrix@=6pt{
	\bar X_1\ar@<-6pt>[dd]_{\bar d_0} \ar@<6pt>[dd]^{\bar d_1} \ar[rrr]^{\check{\lambda}_1}  &&& T(X_1)  \ar@<-10pt>[dd]_{T(d_0)} \ar@<4pt>[drr]^{T(\delta_1)} \\
	&&&&& T^2(X_0) \ar[dll]^{\mu_{X_0}}\\
	X_0 \ar@{ >->}[uu]|{\bar s_0} \ar[rrr]_{\lambda_{X_0}}  
	&&& T(X_0)     \ar@{ >->}[uu]|{T(\bar s_0)}
}
$$
It determines the left hand side $T$-functor in $T$-$Cat\EE$ where $\bar\lambda_1$ is the factorization of $\check\lambda_1$ through $\lambda_{X_1}$:
$$\xymatrix@=6pt{
\bar X_1 \ar[rrr]_{\bar{\lambda}_1}\ar@<2ex>[rrrrrrrr]^{\check{\lambda}_1} \ar@<-6pt>[dd]_{\bar d_0} \ar@<6pt>[dd]^{\bar d_1} &&&	 X_1\ar@<-6pt>[dd]_{d_0}  \ar[rrrrr]_{\lambda_{X_1}} \ar@{.>}@<6pt>[dd]^{"\delta_1"}  &&&&& T(X_1)  \ar@<-10pt>[dd]_{T(d_0)} \ar@<4pt>[drr]^{T(\delta_1)} \\
 &&&	&&&&&&& T^2(X_0) \ar[dll]^{\mu_{X_0}}\\
 X_0\ar@{ >->}[uu]|{\bar s_0} \ar@{=}[rrr] &&&	X_0 \ar@{ >->}[uu]|{s_0} \ar[rrrrr]_{\lambda_{X_0}}  
 	&&&&& T(X_0)     \ar@{ >->}[uu]|{T(\bar s_0)}
}
$$
since $"\delta_1".\bar\lambda_1=\bar d_1$ in $KlT$ means $\delta_1.\bar\lambda_1=\lambda_{X_0}.\bar d_1$ in $\EE$ which is true since $\lambda_{X_0}.\bar d_1=\mu_{X_0}.T(\delta_1).\check\lambda_1=\mu_{X_0}.T(\delta_1).\lambda_{X_1}.\bar\lambda_1=\mu_{X_0}.\lambda_{T(X_0)}.\delta_1.\bar\lambda_1=\delta_1.\bar\lambda_1$.

Now, let $h_{\bullet}: Z_{\bullet} \to X^T_{\bullet}$ be any $T$-functor with $Z_{\bullet} \in Cat\EE$. This means that we get $\delta_1.h_1=\lambda_{X_0}.h_0.d_1$ in $\EE$.
We have to check that the map $h_1:Z_1\to X_1$ factors through $\bar{\lambda}_1$ or equivalently that $\lambda_{X_1}.h_1$ factors through $\check\lambda_1$. For that it is enough to check that we have an internal functor in $Cat\EE$:
$$\xymatrix@=6pt{
	Z_1\ar@<-6pt>[dd]_{d_0} \ar@<6pt>[dd]^{d_1} \ar[rrr]^{\lambda_{X_1}.h_1}  &&& T(X_1)  \ar@<-10pt>[dd]_{T(d_0)} \ar@<4pt>[drr]^{T(\delta_1)} \\
	&&&&& T^2(X_0) \ar[dll]^{\mu_{X_0}}\\
	Z_0 \ar@{ >->}[uu]|{s_0} \ar[rrr]_{\lambda_{X_0}.h_0}  
	&&& T(X_0)     \ar@{ >->}[uu]|{T(\bar s_0)}
}
$$
which is straightforward. The functor $\RR$ preserves the fully faithful functors since the functor $K_T:KlT\to \Alg T$ is cartesian as soon as so is the monad. 
\endproof
Now, when the monad is cartesian, the injective functor:\\ $\TC_{\EE}: \AlgT \into {\rm T}$-Cat$\EE \into Cat(KlT)$ is defined by the following diagram:
$$\xymatrix@=18pt{
	T^2(X)\ar@<-4ex>[rr]_{T(\xi)} \ar[rr]|{\mu_X} \ar@<4ex>@{.>}[rr]^{"1"_{T^2(X)}}  &&	T(X) \ar@<-2ex>[rr]_{\xi} 
	\ar@<2ex>@{.>}[rr]^{"1"_{T(X)}} \ar@<2ex>[ll]|{T(\lambda_X)} \ar@<-2ex>[ll]|{\lambda_{T(X)})} && X\ar[ll]|{\lambda_X}
}
$$
which makes commute the following diagram:
$$\xymatrix@=20pt{
	\AlgT  \ar@{ >->}[r]^{} \ar@<-2pt>[d]_{U^T} & Cat(KlT)  \ar@<2pt>[d]^{(\;)_0} \\
	\EE   \ar@{=}[r]  & \EE \
}
$$
\begin{prop}
	The upper injective functor is cartesian. We get $\RR(\TC_{\EE}(X,\xi))=\Delta.U^T(X,\xi)=\Delta_X$.
\end{prop}
\proof
The first point is straightforward since the pullbacks of internal categories are levelwise and that it is also the case for the $T$-algebras when $T$ is cartesian. For the second point, given any $T$-algebra $(X,\xi)$, consider the following joint pullback diagram:
$$\xymatrix@=5pt{
	X \ar[dd]_{\lambda_X} \ar[rrrr]^{1_X}	&&    &&	X \ar[dd]^{\lambda_X}   \\
	\\
	T(X)	\ar[dd]_{\xi}\ar[rr]_{\lambda_{T(X)}} && T^2(X) \ar[rr]_{\mu_X} \ar[dd]^{T(\xi)} && T(X)    \\
	\\
	X \ar[rr]_{\lambda_X} && T(X)
}
$$ 
\endproof
The composition functor $\AlgT \stackrel{}{\rightarrowtail} Cat(KlT) \stackrel{Cat(\bar F^T)}{\rightarrowtail} Cat(\Alg T)$ gives rise, for any $X\in \EE$ to the following internal category  in $\AlgT$:
$$\xymatrix@=18pt{
	(T^3(X),\mu_{T^2(X)})\ar@<-4ex>[rr]_{T^2(\xi)} \ar[rr]|{T(\mu_X)} \ar@<4ex>[rr]^{\mu_{T(X)}}  &&	(T^2(X),\mu_{T(X)}) \ar@<-2ex>[rr]_{T(\xi)} 
	\ar@<2ex>[rr]^{\mu_X} \ar@<2ex>[ll]|{T^2(\lambda_X)} \ar@<-2ex>[ll]|{T(\lambda_{T(X)})} && (T(X),\mu_X )\ar[ll]|{T(\lambda_X)}
}
$$
For further developments on these internal categories, see \cite{Bat3}. On the other hand, the "intersection" of the inclusions $Cat(\bar F^T): Cat\EE \into Cat(KlT)$ and $\Alg T\into Cat(KlT)$ is clearly the empty set.

\subsection{$T$-groupoids}\label{Tgroup}

Is the notion of $T$-groupoid meaningful? Probably not in general, but it is clear that when a $T$-category coincides with an internal category in $KlT$, it is legitimate to say that a $T$-category:
$$\xymatrix@=18pt{
	X_{\bullet}:&&	X_3 \ar@<-7ex>[rr]_{d_0^2} \ar@<7ex>@{.>}[rr]^{"\delta_3"} \ar@<-3ex>[rr]|{d_1^2} \ar@<3ex>[rr]|{d_2^2}&& X_2\ar@<-4ex>[rr]_{d_0^1} \ar[rr]|{d_1^1} \ar@<4ex>@{.>}[rr]^{"\delta_2"} \ar[ll]|{s_1^2} \ar@<5ex>[ll]|{s_0^2} \ar@<-5ex>[ll]|{s_2^2} &&	X_1 \ar@<-2ex>[rr]_{d_0} 
	\ar@<2ex>@{.>}[rr]^{"\delta_1"} \ar@<2ex>[ll]|{s_0^1} \ar@<-2ex>[ll]|{s_1^1} && X_0 \ar[ll]|{s_0}
}
$$
is a $T$-groupoid when:
$$\xymatrix@=18pt{
	 X_2\ar@<2ex>[rr]^{d_1^1}  \ar@<-2ex>[rr]_{d_0^1}  &&	X_1 \ar[rr]_{d_0} 
	 \ar[ll]|{s_0^1}  && X_0 
}
$$
is a kernel equivalence relation in $\EE$.

\subsection{When $T$-algebras produce $T$-groupoids}

We shall try now to answer the question: when is the image $\TC_{\EE}(X,\xi)$ an internal groupoid in $KlT$?

\begin{prop}\label{ccartalg}
	Let $\Sigma$ be a pullback stable class of morphisms and  $(T,\lambda,\mu)$ be a strongly $\Sigma$-cartesian monad on $\EE$. Suppose that the object $X$ is such that  the folllowing diagram is a kernel equivalence relation:
	$$\xymatrix@=5pt{
	 T^3(X) \ar@<-2ex>[rrrr]_{\mu_{T(X)}} \ar@<2ex>[rrrr]^{T(\mu_{X})} &&&& T^2(X)\ar[rrrr]_{\mu_X}\ar[llll]|{T(\lambda_{T(X)})}&&&& T(X)		
	}
	$$
	Then the object $T(X)$ satisfies the same property. Suppose moreover that 
	$\Sigma$ is a bijective on objects left cancellable subcategory of $\EE$ containing all the isomorphisms.
	Then any algebra $\xi: T(X)\to X$ on $X$  belongs to $\Sigma$ and the following diagram is a kernel equivalence relation:
	$$\xymatrix@=5pt{
	T^2(X) \ar@<-2ex>[rrrr]_{\mu_{X}} \ar@<2ex>[rrrr]^{T(\xi)} &&&& T(X)\ar[rrrr]_{\xi}\ar[llll]|{T(\lambda_{X})}&&&& X		
}
$$
So, the $T$-category $\TC_{\EE}(X,\xi)$ is actually a $T$-groupoid, i.e. an internal groupoid in the Kleisli category $KlT$.
\end{prop}
\proof
Since $\mu_X$ is in $\Sigma$ and $\mu$ is $\Sigma$-cartesian, the map $T(\mu_{T(X)})$ delineates the composition map of an internal category in $\EE$:
$$\xymatrix@=20pt{
	T^4(X)\ar@<-2ex>[rr]_{T^2(\mu_X)} \ar[rr]|{T(\mu_{T(X)})} \ar@<2ex>[rr]^{\mu_{T^2(X)}}  &&	T^3(X) \ar@<-2ex>[rr]_{T(\mu_{X})} 
	\ar@<2ex>[rr]^{\mu_{T(X)}}  && T^2(X)\ar[ll]|{T(\lambda_{T(X)})}
}
$$
When the diagram in question is a kernel equivalence relation, its image by $T$ is a kernel equivalence relation (i.e. the pair $(T(\mu_{T(X)}),T^2(\mu_X))$ is the kernel pair of $T(\mu_X)$), and this category is actually a groupoid. By duality, the pair $(\mu_{T^2(X)},T(\mu_{T(X)}))$ is the kernel equivalence relation of $\mu_{T(X)}$ and we get the first assertion.

Any $T$-algebra $\xi: T(X)\to X$ produces the following diagram:
$$\xymatrix@=18pt{
	T^4(X) \ar[dd]_{T^3(\xi)}\ar@<-2ex>[rr]_{T^2(\mu_X)} \ar[rr]|<<<<<<<{T(\mu_{T(X)})} \ar@<2ex>[rr]^{\mu_{T^2(X)}}  &&	T^3(X) \ar[dd]^{T^2(\xi)}\ar@<-2ex>[rr]_{T(\mu_{X})} 
	\ar@<2ex>[rr]^{\mu_{T(X)}}  && T^2(X)\ar[ll]|>>>>>>>>>{T(\lambda_{T(X)})} \ar[dd]^{T(\xi)}\\
	&&&&\\
	T^3(X)\ar@<-2ex>[rr]_{T^2(\xi)} \ar[rr]|{T(\mu_X)} \ar@<2ex>[rr]^{\mu_{T(X)}}  &&	T^2(X) \ar@<-2ex>[rr]_{T(\xi)} 
	\ar@<2ex>[rr]^{\mu_X}  && T(X)\ar[ll]|>>>>>>>>{T(\lambda_X)}
}
$$
Now, since the square in question is a pullback, the maps $T(\xi)$ and $T^2(\xi)$ makes the lower row an internal groupoid.  So, 
the involutive "inversion" mapping $\gamma_{\xi}: T^2(X)\to T^2(X)$ of this groupoid exchanges the maps $\mu_X$ and $T(\xi)$. Since, by assumption, the isomorphism $\gamma_X$ is in $\Sigma$ which is a subcategory of $\EE$, then $T(\xi)=\mu_X.\gamma_{\xi}$ belongs to $\Sigma$. Now, the identity $\lambda_X.\xi=T(\xi).\lambda_{T(X)}$ shows that $\xi\in\Sigma$, since $\Sigma$ is left cancellable.

It remains to check the last assertion. Since the lower row of the diagram above is a groupoid, by duality, the following diagram is a kernel equivalence relation:
$$\xymatrix@=12pt{
	T^3(X)\ar@<-1ex>[rr]_{T^2(\xi)} \ar@<1ex>[rr]^{T(\mu_X)}   &&	T^2(X) \ar[rr]_{T(\xi)} 
	&& T(X)
}
$$
The endofunctor $T$ is conservative since $\lambda$ is the equalizer of the pair $(\lambda_T,T(\lambda))$ (since $\lambda$ is in $\Sigma$ and $\lambda$ $\Sigma$-cartesian); accordingly it reflects the pullbacks of maps in $\Sigma$, and the following diagram is thus a kernel equivalence relation:
$$\xymatrix@=12pt{
	T^2(X)\ar@<-1ex>[rr]_{T(\xi)} \ar@<1ex>[rr]^{\mu_X}   &&	T(X) \ar[rr]_{\xi} 
	&& X
}
$$
\endproof
\begin{coro}\label{TTTTT_X}
	Let  $(T,\lambda,\mu)$ be a cartesian monad on $\EE$. Suppose the object $X$ is such that the following diagram is a kernel equivalence relation:
	$$\xymatrix@=3pt{
	 T^3(X) \ar@<-2ex>[rrrr]_{\mu_{T(X)}} \ar@<2ex>[rrrr]^{T(\mu_{X})} &&&& T^2(X)\ar[rrrr]_{\mu_X}\ar[llll]|{T(\lambda_X)}&&&& T(X)
	}
	$$
	Then any $T$-algebra $\xi: T(X)\to X$ on $X$ is such that the $T$-category $\TC_{\EE}(X,\xi)$ is actually a $T$-groupoid, i.e. an internal groupoid in $KlT$.
\end{coro}
  In Section \ref{TTTT_X} we shall produce a cartesian monad where this  condition is satisfied for any object $X$.
  
  \section{When internal categories in $\EE$ coincide with $\GG$-categories in $\Pt\EE$}
  
  In this section we shall show that the category $Cat\EE$ of internal categories in $\EE$ coincides with a specific subcategory of the category of $\GG$-categories in $\Pt\EE$.
  
  Recall that the monad $(\GG,\sigma,\pi)$ is strongly $\P$-cartesian 
  on $\Pt\EE$ (Section \ref{groupoid}), accordingly we get:
  \begin{prop}
  	The full subcategory $\GG_{\P}$-$Cat\Pt\EE$ of $\GG$-$Cat\Pt\EE$ whose objects are the $\GG$-categories with a $\P$-cartesian $0$-leg coincides with the category whose objects the discrete fibrations $h_{\bullet}: X_{\bullet}\to Y_{\bullet}$ in $\EE$, where $h_0:X_0\to Y_0$ is endowed with a given splitting $t_0$ and whose morphisms are the commutative squares between discrete fibrations in $Cat\EE$ as on the left hand side:
  	$$\xymatrix@=6pt
  	{ X_{\bullet}\ar[rr]^{\psi_{\bullet}}  \ar[dd]_{h_{\bullet}}&&
  		\bar X_{\bullet}    \ar[dd]^{\bar h_{\bullet}} &&&& X_{0}\ar[rr]^{\psi_{0}}  \ar@<-1ex>[dd]_{h_{0}}&&
  		\bar X_{0}    \ar@<-1ex>[dd]_{\bar h_{0}} \\
  		&&&& &&\\
  		Y_{\bullet}  \ar[rr]_{\phi_{\bullet}} &&
  		\bar Y_{\bullet}  &&&&   Y_{0}  \ar[rr]_{\phi_{0}}\ar[uu]_{t_0} &&
  		\bar Y_{0} \ar[uu]_{\bar t_0}
  	}
  	$$
  	such that the above right hand side square is a morphism in $\Pt\EE$.
  \end{prop}
  \proof
  Apply Propositions \ref{ssigma} and \ref{KlG}.
  \endproof
  We shall be interested now in the full subcategory of  $\GG_{\P}$-$Cat\Pt\EE$ whose objects are the $\GG$-categories having, in addition, an idomorphic $1$-leg.
  \begin{theo}\label{mmmain}
  	The full subcategory of  $\GG$-$Cat\Pt\EE$ whose objects are the $\GG$-categories with a $\P$-cartesian $0$-leg and an idomorphic $1$-leg is isomorphic to the category $Cat\EE$ of internal categories in $\EE$.
  \end{theo} 
  \proof
  Let $h_{\bullet}: X_{\bullet}\to Y_{\bullet}$ be a $\GG$-category with a $\P$-cartesian $0$-leg, namely a vertical discrete fibration with a section $t_0$ of $h_0$:
  $$\xymatrix@=12pt
  { X_0 \ar[rr]|{s_0^Y}  \ar@<-1ex>[dd]_{h_0}&&
  	X_1    \ar[dd]^{h_1} \ar@<2ex>[ll]^{d_{1^X}} \ar@<-2ex>[ll]_{d_{0^X}} && X_2 \ar[dd]^{h_2}  \ar@<2ex>[ll]^{d_{2}}\ar@<-2ex>[ll]_{d_{0}}\ar[ll]|{d_{1}}\\
  	&& &&\\
  	Y_0  \ar[rr]|{s_0} \ar[uu]_{t_0} &&
  	Y_1     \ar@<2ex>[ll]^{d_{1^Y}} \ar@<-2ex>[ll]_{d_{0^Y}} && Y_2  \ar@<2ex>[ll]^{d_{2}}\ar@<-2ex>[ll]_{d_{0}}\ar[ll]|{d_{1}}
  }
  $$
  we shall denote by $t_i$ the induced section of $h_i$.
  Its underlying $\GG$-graph: $$(h_0,t_0) \stackrel{(d_0^Y,d_0^X)}{\longleftarrow} (h_1,t_1) \stackrel{(d_1^X,R(d_1^X)).\sigma(h_1,t_1)}{\longrightarrow} \GG(h_0,t_0)$$ is the following one:
  $$\xymatrix@=8pt
  { 
  	&& 	X_1    \ar@<-1ex>[dd]_{h_1}  \ar[lld]_{d_{0^X}} \ar[rrd]^{(d_1^X.t_1.h_1,d^X_1)}\\
  	X_0   \ar@<-1ex>[dd]_{h_0}&&
  	&& R[h_0] \ar@<-1ex>[dd]_{p_0^{h_0}}  \\
  	&& Y_1  \ar[lld]^{d_{0^Y}}\ar[uu]_{t_1} \ar[rrd]_{d_1^X.t_1}\\ 
  	Y_0   \ar[uu]_{t_0} &&
  	&& X_0 \ar[uu]_{s_0^{h_0}} 
  }
  $$
  Saying that its $1$-leg is idomorphic is saying that $d_1^X.t_1=1_{Y_1}$. Whence:\\
  1) $h_0=h_0.d_1^X.t_1=d_1^Y.h_1.t_1=d_1^Y$, and 2) $t_0=t_0.d_1^Y.s_0^Y=d_1^X.t_1.s_0^Y=s_0^Y$, and consequently $(h_0,t_0)=(d_1,s_0)$. Similarly we have $d_2^X.t_2=1_{Y_2}$. Whence:\\
  1) $h_1=h_1.d_2^X.t_2=d_2^Y.h_2.t_2=d_2^Y$, and 2) $t_1=t_1.d_2^Y.s_1^Y=d_2^X.t_2.s_1^Y=s_1^Y$, and consequently $(h_1,t_1)=(d_2,s_1)$. Accordingly, we get the following diagram:
  $$\xymatrix@=12pt
  { Y_1 \ar[rr]|{s_0}  \ar@<-1ex>[dd]_{d_1}&&
  	Y_2    \ar@<-1ex>[dd]_{d_2} \ar@<2ex>[ll]^{d_{1}} \ar@<-2ex>[ll]_{d_{0}} && Y_3 \ar@<-1ex>[dd]_{d_3}  \ar@<2ex>[ll]^{d_{2}}\ar@<-2ex>[ll]_{d_{0}}\ar[ll]|{d_{1}}\\
  	&& &&\\
  	Y_0  \ar[rr]|{s_0} \ar[uu]_{s_0} &&
  	Y_1     \ar@<2ex>[ll]^{d_{1}} \ar@<-2ex>[ll]_{d_{0}}  && Y_2  \ar@<2ex>[ll]^{d_{2}}\ar@<-2ex>[ll]_{d_{0}}\ar[ll]|{d_{1}}
  }
  $$
  which is nothing but the discrete fibration $\epsilon_{Y_{\bullet}}: DecY_{\bullet} \to Y_{\bullet}$ with the section $s_0$, and nothing more. So that a $\GG_\P$-category with an idomorphic $1$-leg is just an internal category in $\EE$. Conversely any internal category $Y_{\bullet}$ in $\EE$ produce the above $\GG_\P$-category with an idomorphic $1$-leg.
  \endproof
  Now, consider the injective functor $\GGC: \Alg\GG=Grd\EE \into \GG$-$Cat(\Pt\EE)$; according to Proposition \ref{ccartalg}, any $\GG$-algebra is $\P$-cartesian, so that the pointed $\GG$-graph underlying $\GGC((d_0,s_0),(d_1,d_2))$ (following the notations of Section \ref{groupoid}):
  $$\xymatrix@=10pt{
  	& \GG(d_0,s_0) \ar[dl]_{(d_1,d_2)}\ar[dr]^{1_{\GG(d_0,s_0)}}\\
  	(d_0,s_0) \ar[rr]_{\sigma(d_0,s_0)} && \GG(d_0,s_0)
  }
  $$
  has a $\P$-cartesian $0$-leg and an idomorphic $1$-leg.
  Accordingly, the injective functor $\GGC$ factors through $Cat\EE$, producing the natural inclusion $Grd\EE\into Cat\EE$. So, not only the monad  $(\GG,\sigma,\pi)$ produces the category $Grd\EE=\Alg\GG$ of internal groupoids, but also it entirely rules the previous inclusion.
  
  \section{Internal $n$-groupoids and $n$-categories}
  
  In this section, we shall show that the constructions and results of the previous section about the monad $(\GG,\sigma,\pi)$ have a natural extension to the internal (strict) $n$-groupoids and $n$-categories.
  
  \subsection{The monad $(\GG_F,\sigma_F,\pi_F)$}\label{locfib}
  
  We shall first introduce a locating process for the monad $(\GG,\sigma,\pi)$
   up to a fibration. So let $F:\bar\EE\to\EE$ be any fibration whose underlying functor is cartesian; we denote by $F_W$ the fiber above $W\in\EE$. 
   \begin{lemma}
   Given any fibration $F:\bar\EE\to\EE$ whose underlying functor is cartesian, then any fiber $F_W$ has pullbacks. 
   \end{lemma}
\proof Given pair $(f,g)$ of maps with same codomain in the fiber $F_{F(Y)}$ and their pullback in $\bar\EE$:
    $$\xymatrix@=6pt
  { P \ar[rr]^{p_X}  \ar[dd]_{p_Y}&&
  	X \ar[dd]^{f} &&&&&  F(P) \ar@<-1ex>[dd]_{F(p_Y)}  \ar[rr]^{F(p_X)} && F(X) \ar@{=}[dd]^{F(f)}  \\
  	&& &&&&&&\\
  	Z  \ar[rr]_{g}  &&
  	Y  &&&&&   F(Z) \ar@{=}[rr]_{F(g)}  && F(Y)  
  }
  $$
  we get $F(p_X)=F(p_Y)=\gamma$, where $\gamma$ is an isomorphism in $\EE$. Taking $\zeta: \bar P\to P$ the cartesian map  above $\gamma^{-1}$ with domain $P$ furnishes the pullback of the pair $(f,g)$ inside the fiber $F_{F(Y)}$. 
  \endproof
  Let us denote by $\Pt_F\bar\EE$ the full subcategory of $\Pt\bar\EE$ whose objects are the split epimorphisms in a fiber of $F$, it is obtained by the following left hand side pullback where $1_\EE(W)=(1_W,1_W)$:
  $$\xymatrix@=8pt
  { \Pt_F\bar\EE \ar@{ >->}[rr]_{\iota_{\bar\EE}} \ar@<2ex>[rrrrr]^{\P_F}  \ar[dd]_{} &&
  	\Pt\bar\EE \ar[dd]^{\Pt F} \ar[rrr]_{\P_{\bar\EE}} &&& \bar\EE \\
  	&& &&&&&&\\
  	\EE  \ar@{ >->}[rr]_{1_\EE}  && \Pt\EE   
  }
  $$
  This makes $\Pt_F\bar\EE$ a cartesian category and any functor in this diagram is a cartesian one. This produces, in addition, the upper horizontal cartesian functor $\P_F$ which becomes a subfibration of $\P_{\bar\EE}$ by the following \emph{specification of the base-change functor $\psi^*$}: starting with any split epimorphism $(f,s):X\splito Y$ in the fiber $F_{F(Y)}$, take the pullback in $\bar\EE$ as on the left hand side which produces a cartesian map above $\psi$ in $\Pt\bar\EE$:
   $$\xymatrix@=6pt
  { \bar X \ar[rrr]^{\bar\psi}  \ar[dd]_{\bar f}&&&
  	X \ar[dd]_{f} &&&&&  F(\bar X) \ar@<-1ex>[dd]_{F(\bar f)}  \ar[rr]^{F(\bar\psi)} && F(X) \ar@{=}[dd]^{F(f)}  \\
  	&& &&&&&&\\
  	Z  \ar[rrr]_{\psi} \ar@<-1ex>[uu]_{\bar s} &&&
  	Y \ar@<-1ex>[uu]_s &&&&&   F(Z) \ar[rr]_{F(\psi)}  && F(Y)  
  }
  $$
  So, in the right hand side image by $F$ in $\EE$ which is a pullback, the map $F(\bar f)=\gamma$ is an isomorphism whose inverse is $F(\bar s)$. Taking the cartesian isomorphism $\zeta:\check X \to \bar X$ above $\gamma^{-1}$ with codomain  $\bar X$ produces the desired split epimorphism $(\bar f.\zeta, \zeta^{-1}.\bar s): \check X\splito Z$ in the fiber above $F(Z)$. \emph{From now on we shall use the previous specification in the construction of the base-change functors $\psi^*$}. 
  
  We shall denote by  $\P_F$ the class of the cartesian maps with respect to the fibration $\P_F$ (namely pullbacks between split epimorphisms belonging to a fiber) and, again, we shall call idomorphims the morphisms $(y,x)$ in $\Pt_F\bar\EE$ whose lower map $y$ in $\bar\EE$ is an identity map. Modulo the above precisions, \emph{the monad $(\GG,\sigma,\pi)$ on $\Pt{\bar\EE}$ is stable on $\Pt_F\bar\EE$}; for sake of clarity, we shall denote it by $(\GG_F,\sigma_F,\pi_F)$.
  \begin{prop}\label{richF}
  	The endofunctor $\GG_F$ on $\Pt_F\bar\EE$ is cartesian. It preserves and reflects the class $\P_F$. The monad $(\GG_F,\sigma_F,\pi_F)$ is strongly $\P_F$-cartesian. Furthermore, given any object $(g,t)$ in $\Pt_F{\bar\EE}$, the following diagram is a kernel equivalence relation in $\Pt\bar\EE$ with its (levelwise) quotient:
  	$$\xymatrix@=25pt
  	{
  		\GG_F^3(g,t) \ar@<-2ex>[rr]_{\GG_{F\pi_F{(g,t)}}} \ar@<3ex>[rr]^{\pi_F{\GG_F(g,t)}} && \GG_F^2(g,t) \ar[r]_{\pi_{F(g,t)}} \ar[ll]_{\GG_F(\sigma_{F\GG_F(g,t)})}  & \GG_F(g,t) 
  	}
  	$$
  \end{prop}
\proof
It is just Proposition \ref{rich} restricted to the full subcategory  $\Pt_F{\bar\EE}$ of $\Pt\bar\EE$ since the inclusion $\iota_{\bar\EE}$ preserves the cartesian maps and the monad $(\GG_F,\sigma_F,\pi_F)$ is just the restriction to $\Pt_F{\bar\EE}$ of the monad $(\GG,\sigma,\pi)$ on $\Pt\bar\EE$.
\endproof
\begin{prop}\label{caracF}
	Any algebra $\alpha:\GG_F(g,t)\rightarrow (g,t)$ of this monad necessarily  belongs to the class $\P_F$. The category of algebras of the monad $(\GG_F,\sigma_F,\pi_F)$ on $\Pt_F\bar\EE$ is the full subcategory $Grd^F\bar\EE$ of $Grd\bar\EE$ whose objects are the internal groupoids  in the fibers of $F$.
\end{prop}
\proof
This time, it is just a restriction of Theorem \ref{carac} to the full subcategory $\Pt_F\bar\EE$ of $\Pt\bar\EE$.
\endproof
We shall denote by $(\;)^F_0$ the diagonal functor of the following commutative square:
$$\xymatrix@=6pt
{ \Alg \GG_F=Grd^F\bar\EE \ar@{ >->}[rrrr]^{} \ar@{-->}[rrrrddd]^>>>>>>>>>>>>{(\;)^F_0}  \ar[ddd]_{U^{\GG_F}} &&&&
	Grd\bar\EE \ar[ddd]^{(\;)_0} \\
	&& &&&&&&\\
	&& &&&&&&\\
	\ar@{.>}@(ul,dl)_{\GG_F}\Pt_F\bar\EE  \ar[rrrr]_{\P_F}  &&&& \bar\EE  
}
$$
\begin{defi}
	A functor $F:\bar\EE \to \EE$ is called a fibered reflection when it is cartesian and is a fibration such that any fiber $F_W$ of $F$ has a terminal object $T(W)$ which is stable under any base-change functor.
\end{defi}
The easiest examples of fibered reflection are the fibrations $(\;)_0: Cat\EE \to \EE$ and $(\;)_0: Grd\EE \to \EE$, when $\EE$ is a cartesian category with a terminal object; in both cases, the  terminal object in the fiber above the object $X$ in $\EE$ being the undiscrete equivalence relation $\nabla_X=R[\tau_X]$, where $\tau_X: X\to 1$ is the terminal map. We are going to show now that when $F$ is a fibered reflection, so are the forgetful functors $(\;)^F_0: Grd^F\EE\to \EE$ and $(\;)^F_0: Cat^F\EE\to \EE$. First, the previous terminology comes from the following:
\begin{prop}\cite{B2}
	A functor $F:\bar\EE \to \EE$ is a fibered reflection if and only if:\\
	1) the functor $F$ is cartesian and has a right adjoint right inverse $T:\EE\to\bar\EE$ such that the unit $\eta_X: X\to TF(X)$ of this co-adjoint pair is such that $F(\eta_X)=1_{F(X)}$;\\
	2) for any map $h: Z\to F(X)$ in $\bar\EE$, there is a map $\bar h: \bar Z\to X$ in $\EE$ such that $F(\bar h)=h$ and the following square is a pullback:
	$$\xymatrix@=6pt
	{ \bar Z \ar[rrrr]^{\bar h}  \ar[dd]_{\eta_{\bar Z}} &&&&
		X \ar[dd]^{\eta_X} \\
		&& &&&&&&\\
		TF(\bar Z)=T(Z)  \ar[rrrr]_>>>>>>>>{TF(\bar h)=T(h)}  &&&& TF(X)   
	}
	$$ 
\end{prop}
\proof
Suppose $F$ is a fibered reflection. Choose a terminal object $T(W)$ in the fiber $F_W$; this determines a right adjoint right inverse $T$ of $F$ with $\eta_X: X\to TF(X)$ the terminal map in the fiber $F_{F(X)}$. Then a map $f: X\to Y$ in $\bar\EE$ is cartesian with respect to $F$ if and only if the following left hand side square is a pullback in $\bar\EE$, see Section 1 in \cite{B2} for instance:
$$\xymatrix@=6pt
{ X \ar[rrr]^{f}  \ar[dd]_{\eta_X} &&&
	Y \ar[dd]^{\eta_Y} &&&  \check Z \ar[rrr]^{\check h}  \ar[dd]_{\check\eta} &&&
	X \ar[dd]^{\eta_X}   \\
	&& &&&&&&\\
	TF(X)  \ar[rrr]_{TF(f)}  &&& TF(Y) &&& 	T(Z)  \ar[rrr]_{T(h)}  &&& TF(X)   
}
$$
Now starting with a map $h: Z\to F(X)$ in $\bar\EE$, take the above right hand side pullback in $\EE$. Since $F$ is cartesian, its image by $F$ is a pullback in $\EE$, and $F(\check\eta)$ is an isomorphism. Taking $\gamma: \bar Z\to\check Z$ the invertible cartesian map above $F(\check\eta)^{-1}$ with codomain $\check Z$,
produces the desired map $\bar h=\check h.\gamma:\bar Z\to X$ of condition 2).

Conversely, let $F$ be a functor satisfying the two above conditions. Condition 1) implies that any map $f: X\to Y$ making the above left hand side square a pullback is cartesian with respect to $F$, while Condition 2) guarantees the existence of a cartesian map above any map $h$. Then $\eta_X: X\to TF(X)$ is necessarily the terminal map in the fiber $F_{F(X)}$.
\endproof
\begin{prop}\label{fiberef}
	When the fibration $F:\bar\EE \to \EE$ is a fibered reflection, so is the forgetful functor $(\;)^F_0: Grd^F\bar\EE\to \bar\EE$.
\end{prop}
\proof
The kernel equivalence relation $R[\eta_X]$ produces a groupoid $R_{\bullet}[\eta_X]$ in the fiber $F_{F(X)}$ (since $F(\eta_X)=1_{F(X)}$) such that $R_{0}[\eta_X]=X$, and it is clearly a terminal object among the groupoids $Z_{\bullet}$ in the fiber $F_{F(X)}$ such that $Z_0=X$. 

Let us check that $(\;)^F_0$ is a fibration whose base-change functors preserve these terminal objects. So, let $X_{\bullet}$ be any internal groupoid in a fiber of $F$ and $h:Z\to X_0$ be any map in $\bar\EE$; then consider the following left hand side pullback in the following left hand side diagram in $\bar\EE$, where $\bar\eta_1X_{\bullet}$ is the factorization of the pair $(d_0^{X_{\bullet}},d_1^{X_{\bullet}})$:
$$\xymatrix@=20pt
{
  \bar Z_1 \ar[r]^{\bar\eta} \ar[d]_{\bar h_1} &	R[\eta_Z] \ar[d]_-{R(h)} \ar@<1ex>[r]^{p_0} \ar@<-1ex>[r]_{p_1} & Z  \ar[d]^-{h} \ar[r]^{\eta_Z} \ar[l]|{} &
	TF(Z) \ar[d]^-{TF(h)} && \bar Z_{\bullet} \ar[r]^{\bar\eta_{\bullet}} \ar[d]_{\bar h_{\bullet}} & R_{\bullet}[\eta_{Z}] \ar[d]^-{R_{\bullet}(h)}\\
 X_1 \ar[r]_{\bar\eta_1X_{\bullet}}	& R[\eta_{X_0}] \ar@<1ex>[r]^{p_0} \ar@<-1ex>[r]_{p_1} & X_0 \ar[r]_-{\eta_{X_0}} \ar[l] & TF(X_0) && X_{\bullet} \ar[r]_{\bar\eta_{\bullet}X_{\bullet}} & R_{\bullet}[\eta_{X_0}] 
}
$$
This produces an internal groupoid in $\bar\EE$, since this pullback in $\bar\EE$ is underlying the right hand side pullback in $Grd\bar\EE$. The map $\bar\eta$ is not necessarily inside a fiber, but certainly $F(\bar\eta)$ is an isomorphism since $F(\bar\eta_1X_{\bullet})=1_{F(X_0)}$. Take the invertible cartesian map $\zeta_1: Z_1\to \bar Z_1$ above  $F(\bar\eta)^{-1}$ with codomain $\bar Z_1$, then the associated internal groupoid $Z_{\bullet}$ (which is isomorphic to $\bar Z_{\bullet}$) belongs to $Grd_F\bar\EE$. It then straighforward to check that the internal functor $h_{\bullet}.\zeta_{\bullet}: Z_{\bullet}\to X_{\bullet}$ is the desired cartesian map above $h$ with respect to the functor: $(\;)^F_0: Grd^F\bar\EE\to \bar\EE$. This construction makes $R_{\bullet}(h)$ a cartesian map above $h$, which means that the terminal object $R_{\bullet}[\eta_{X_0}]$ in the fiber of $(\;)^F_0$ is preserved by the base-change functor along $h$; in other words this means that the fibration $(\;)^F_0$ is a fibered reflection.
\endproof
\begin{theo}\label{mmmainF}
	The full subcategory of $\GG_F$-$Cat(\Pt_F\bar\EE)$ whose objects are the $\GG_F$-categories with a $\P_F$-cartesian $0$-leg  and an idomorphic $1$-leg is isomorphic to the full subcategory $Cat^F\bar\EE$ of $Cat\bar\EE$ whose objects are the internal categories in the fibers of $F$. The inclusion $\GGC_F: \Alg\GG_F \into Cat^F\bar\EE$ coincides with the following upper one:
	$$\xymatrix@=20pt
	{
		Grd^F\bar\EE \ar@{ >->}[r] \ar[d]_{(\;)^F_0} &	Cat^F\bar\EE \ar[d]^-{(\;)^F_0}  \\
	 \bar\EE \ar@{=}[r] & \bar\EE 
	}
	$$
	When $F$ is a fibered reflection, so is $(\;)^F_0: Cat^F\bar\EE \to\bar\EE$.
\end{theo}
\proof
Once again, the first point is only a straighforward restriction of Theorem \ref{mmain} to the full subcategory $\Pt_F\bar\EE$ of $\Pt\bar\EE$. And the second one holds since the previous proof for $Grd^F\bar\EE$ is still valid for $Cat^F\bar\EE$, which means, as expected, that the terminal objects in the fibers $Grd^F_X\bar\EE$ and $Cat^F_X\bar\EE$ are the same one.
\endproof

\subsection{When the fibration $F$ has protomodular fibers}\label{protfib}

Recall from \cite{B0}, that \emph{a category $\CC$ is protomodular} when any base-change functor of the fibration $\P_\CC: \Pt\CC\to\CC$ is conservative,  that any protomodular category is a Mal'tsev one in the sense of \cite{CLP} (namely any reflexive relation in $\CC$ is an equivalence relation), and that any internal category in a Mal'tsev category is an internal groupoid \cite{CPP}. 

The easiest examples of protomodular category are the category $\Gp$ of groups and $\Gp\EE$ of internal groups in $\EE$ when $\EE$ is cartesian. More generally any fiber of the fibration $(\;)_0:Grd\EE\to\EE$ is protomodular, again see \cite{B0}, and when $\CC$ is protomodular, so is the category $Grd\CC$. 
When the fibration $F$ has protomodular fibers, then the extension determined by Theorem \ref{mmmainF} does not produce anything new:
\begin{prop}\label{mainprot}
	 When the fibration $F$ has protomodular fibers, the inclusion $\GGC_F: \Alg\GG_F=Grd^F\bar\EE \into Cat^F\bar\EE$ is an isomorphism of categories, and the fibered reflection $(\;)^F_0: Grd^F\bar\EE\to \bar\EE$ has protomodular fibers as well.
\end{prop}
\proof
Following what we just recalled, internal categories and internal groupoids  do coincide inside the protomodular fibers of $F: \bar\EE\to\EE$; whence the first point. Now, the fiber $(\;)^{F,X}_0$ of $(\;)^F_0$ above the object $X\in\bar\EE$ is a cartesian subcategory of the category $Grd(F_{F(X)})$ which, as we just recalled above, is protomodular since so is $F_{F(X)}$. Accordingly, so is this fiber $(\;)^{F,X}_0$.
\endproof

\subsection{$2$-categories and $2$-groupoids}

(Strict) $2$-categories and $2$-groupoids have been introduced by Ehresmann in \cite{Eh}. Internally speaking, the category $2$-$Cat\EE$ of internal $2$-categories in $\EE$ is nothing but the full subcategory of the category $Cat(Cat\EE)$ of double categories whose objects are the internal categories in the fibers of the fibration $(\;)_0: Cat\EE\to\EE$, see Section VI.2 in \cite{B2}. So we are in the situation investigated in Section \ref{locfib} with the fibration $F=(\;)_0$. This section will be devoted to the translation  of the results of Section \ref{locfib}, and this will show how, again, the monad $(\GG,\sigma,\pi)$ entirely rules the construction of the category $2$-$Cat\EE$. Similarly, the category $2$-$Grd\EE$ of internal $2$-groupoids in $\EE$ is nothing but the full subcategory of the category $Grd(Grd\EE)$ of double groupoids whose objects are the internal groupoids in the fibers of the fibration $(\;)_0: Grd\EE\to\EE$.

Let us begin by the category $2$-$Grd\EE$. In this way, $2$-$Grd\EE=Grd^{(\;)_0}Grd\EE$. When there is no ambiguity, a $2$-groupoid will be denoted by the central part $X^2_{\bullet}$ of the internal groupoid defining it in a fiber of $(\;)_0: Grd_EE\to\EE$:
$$\xymatrix@=25pt
{
	X^2_{\bullet}\times_0 X^2_{\bullet} \ar[rr]|{d^{2,1}_{\bullet}} \ar@<2ex>[rr]^{d^{2,2}_{\bullet}}\ar@<-2ex>[rr]_{d^{2,0}_{\bullet}} && X^2_{\bullet} \ar@<-2ex>[rr]_{d^{2,0}_{\bullet}} \ar@<2ex>[rr]^{d^{2,1}_{\bullet}} && X^1_{\bullet} \ar[ll]|{s^{2,0}_{\bullet}}  
}
$$
where the left hand side object is a pullback in this fiber. The internal groupoid $X^1_{\bullet}$ is called the groupoid of $1$-morphisms, while the internal groupoid $X^2_{\bullet}$ is called the groupoid of $2$-morphisms or $2$-cells. Let us translate now the results of Section \ref{locfib} with $F=(\;)_0: Grd\EE\to\EE$. For that and for sake of simplicity we shall denote:\\
1) by $\Pt_0Grd\EE$ the category  $\Pt_{(\;)_0}Grd\EE$ whose objects are the split epimorphisms between internal groupoids lying in a fiber of $(\;)_0$,\\ 
2) by $\P_0\EE$ the fibration $\P_{(\;)_0}: \Pt_0Grd\EE=\Pt_{(\;)_0}Grd\EE\to Grd\EE$ associating with any split epimorphism of this kind its codomain,\\
3) and by $(\GG_1,\sigma_1,\pi_1)$ the monad $(\GG_{(\;)_0},\sigma_{(\;)_0},\pi_{(\;)_0})$ on the category $\Pt_0Grd\EE$.\\
So we get:
\begin{prop}\label{carac2g}
	1) The category $2$-$Grd\EE$ is isomorphic to $\Alg \GG_1$.
	 The forgetful funtor $(\;)_1: 2$-$Grd\EE\to Grd\EE$ associating the groupoid $X^1_{\bullet}$ with the $2$-groupoid $X^2_{\bullet}$ is a fibered reflection.\\
	2) The inclusion $\GGC_1: \Alg\GG_1=2$-$Grd\EE \into Cat^{(\;)_0}Grd\EE$ is an isomorphism of categories and the fibered reflection $(\;)_1: 2$-$Grd\EE\to Grd\EE$ has protomodular fibers.
\end{prop}
\proof
For the first point, just apply Propositions \ref{richF}, \ref{caracF} and \ref{fiberef} to the fibration $(\;)_0: Grd\EE \to\EE$. For the second one, apply Proposition \ref{mainprot}.
\endproof

Let us translate now the results of Section \ref{locfib} related with $F=(\;)_0:Cat\EE\to\EE$. In this way, $2$-$Cat\EE=Cat^{(\;)_0}Cat\EE$. Again, when there is no ambiguity, a $2$-category will be denoted by the central part $X^2_{\bullet}$  of the internal category producing it in a fiber of  of $(\;)_0: Cat\EE\to\EE$:
$$\xymatrix@=25pt
{
	X^2_{\bullet}\times_0 X^2_{\bullet} \ar[rr]|{d^{2,1}_{\bullet}} \ar@<2ex>[rr]^{d^{2,2}_{\bullet}}\ar@<-2ex>[rr]_{d^{2,0}_{\bullet}} && X^2_{\bullet} \ar@<-2ex>[rr]_{d^{2,0}_{\bullet}} \ar@<2ex>[rr]^{d^{2,1}_{\bullet}} && X^1_{\bullet} \ar[ll]|{s^{2,0}_{\bullet}}  
}
$$
where the left hand side object is a pullback in this fiber of $(\;)_0: Cat\EE\to\EE$. The internal category $	X^1_{\bullet}$ is called the category of $1$-morphisms, while the internal category $X^2_{\bullet}$ is called the category of $2$-morphisms or $2$-cells.  Again for sake of simplicity, we shall denote:\\
1) by $\Pt_0Cat\EE$ the category  $\Pt_{(\;)_0}Cat\EE$ whose objects are the split epimorphisms between internal categories lying in a fiber of $(\;)_0$,\\
2) by $\P^C_0\EE$ the fibration $\P_{(\;)_0}: \Pt_0Cat\EE=\Pt_{(\;)_0}Cat\EE\to Cat\EE$ associating with any split epimorphism of this kind its codomain,\\
3) and by $(\GG^C_1,\sigma^C_1,\pi^C_1)$ the monad $(\GG_{(\;)_0},\sigma_{(\;)_0},\pi_{(\;)_0})$ on the category $\Pt_0Cat\EE$.\\ Now translating the results of Section \ref{locfib} we get:
\begin{theo}\label{mmmainc1}
	1) The category $\Alg \GG^C_1$ is the full subcategory $2_G$-$Cat\EE$ of $Grd(Cat\EE)$ whose objects are the internal groupoids  in the fibers of $(\;)_0: Cat\EE\to\EE$, namely the $2$-categories with invertible $2$-cells. The forgetful funtor $(\;)_1: 2_G$-$Cat\EE\to Cat\EE$ associating the category $X^1_{\bullet}$ with the $2$-category $X^2_{\bullet}$ is a fibered reflection.\\ 
	2) The full subcategory of $\GG^C_1$-$Cat(\Pt_0Cat\EE)$ whose objects are the $\GG^C_1$-categories  with a $\P^C_0\EE$-cartesian $0$-leg and an idomorphic $1$-leg is isomorphic to the category of internal categories in the fibers of $(\;)_0$, namely to the category $2$-$Cat\EE$ of internal $2$-categories.\\
	3) The inclusion $\GGC_1: \Alg\GG^C_1 \into 2$-$Cat\EE$ coincides with the following upper horizontal one in the following commutative diagram:
	$$\xymatrix@=20pt
	{
		2_{G}{\rm -}Cat\EE \ar@{ >->}[r] \ar[d]_{(\;)_1} &	2{\rm -}Cat\EE \ar[d]^-{(\;)_1}  \\
		Cat\EE \ar@{=}[r] & Cat\EE 
	}
	$$
	where the vertical functors are fibered reflections.
\end{theo}
\proof
Apply Proposition \ref{caracF} and Theorem \ref{mmmainF}.
\endproof

  \subsection{$n$-categories and $n$-groupoids}
 
Again, (strict) $n$-categories have been introduced in \cite{Eh}. Internally speaking, the category $(n+1)$-$Cat\EE$ of internal $(n+1)$-categories in $\EE$ is defined by induction, see for instance \cite{Banother}. 

We defined the fibered reflection $(\;)_1: 2$-$Cat\EE\to Cat\EE$. Suppose we have defined the fibered reflection $(\;)_{n-1}: n$-$Cat\EE\to (n$-$1)$-$Cat\EE$.
Then the category $(n+1)$-$Cat\EE$ of $(n+1)$-categories is the full subcategory of $Cat(n$-$Cat\EE)$ whose objects are the internal categories in the fibers of $(\;)_{n-1}$. We are now in the situation investigated in Section \ref{locfib} with $F=(\;)_{n-1}$. This section will be devoted to the translation of the results of Section \ref{locfib}, and this will show how, this time, the monad $(\GG,\sigma,\pi)$ entirely rules the construction of the category $(n+1)$-$Cat\EE$. Similarly, the category $(n+1)$-$Grd\EE$ of internal $(n+1)$-groupoids in $\EE$ is inductively defined as the full subcategory of the category $Grd(n$-$Grd\EE)$ whose objects are the internal groupoids in the fibers of the fibration $(\;)_{n-1}: n$-$Grd\EE\to (n-1)$-$Grd\EE$.

Let us begin by the category $(n+1)$-$Grd\EE$. When there is no ambiguity, a $(n+1)$-groupoid will be denoted by the central part $X^{n+1}_{\bullet}$ of the internal groupoid defining it in a fiber of $(\;)_{n-1}: n$-$Grd\EE\to (n-1)$-$Grd\EE$:
$$\xymatrix@=25pt
{
	X^{n+1}_{\bullet}\times_{n_1} X^{n+1}_{\bullet} \ar[rr]|{d^{n+1,1}_{\bullet}} \ar@<2ex>[rr]^{d^{n+1,2}_{\bullet}}\ar@<-2ex>[rr]_{d^{n+1,0}_{\bullet}} && X^{n+1}_{\bullet} \ar@<-2ex>[rr]_{d^{n+1,0}_{\bullet}} \ar@<2ex>[rr]^{d^{n+1,1}_{\bullet}} && X^n_{\bullet} \ar[ll]|{s^{n+1,0}_{\bullet}}  
}
$$
where the left hand side object is a pullback in this fiber. The $n$-groupoid $X^n_{\bullet}$ is called the $n$-groupoid of $n$-morphisms, while the $n$-groupoid $X^{n+1}_{\bullet}$ is called the $n$-groupoid of $(n+1)$-morphisms or $(n+1)$-cells. Let us translate now the results of Section \ref{locfib} with $F=(\;)_{n-1}$. For that and for sake of simplicity we shall denote:\\
1) by $\Pt_{n-1}n$-$Grd\EE$ the category  $\Pt_{(\;)_{n-1}}n$-$Grd\EE$ whose objects are the split epimorphisms between internal $n$-groupoids  lying in a fiber of $(\;)_{n-1}$,\\ 
2) by $\P_{n-1}\EE$ the fibration $\P_{(\;)_{n-1}}: \Pt_{n-1}n$-$Grd\EE\to n$-$Grd\EE$ associating with any split epimorphism of this kind its codomain,\\
3) and by $(\GG_n,\sigma_n,\pi_n)$ the monad $(\GG_{(\;)_{n-1}},\sigma_{(\;)_{n-1}},\pi_{(\;)_{n-1}})$ on the category $\Pt_{n-1}n$-$Grd\EE$.
Now translating the results of the Sections \ref{locfib} and \ref{protfib}, we get:
\begin{prop}\label{caracgn}
	1) The category $(n+1)$-$Grd\EE$ of internal $(n+1)$-groupoids is isomorphic to $\Alg \GG_n$. The forgetful functor $(\;)_{n}: (n+1)$-$Grd\EE\to n$-$Grd\EE$ associating the $n$-groupoid $X^n_{\bullet}$ with the $(n+1)$-groupoid $X^{n+1}_{\bullet}$ is a fibered reflection.\\
	2) The inclusion: $\GGC_n: \Alg\GG_{n}=(n+1){\rm -}Grd\EE \into Cat^{(\;)_{n-1}}n{\rm -}Grd\EE$ is an isomorphism of categories and the fibered reflection: $(\;)_{n}: (n+1)$-$Grd\EE\to n$-$Grd\EE$ has protomodular fibers.
\end{prop}
\proof
For the first point, apply Propositions \ref{richF}, \ref{caracF} and \ref{fiberef} to the fibered reflection $(\;)_{n-1}: n$-$Grd\EE \to (n$-$1)$-$Grd\EE$.
For the second one, apply Proposition \ref{mainprot}.
\endproof 

Let us translate now the results of Section \ref{locfib} related to the fibration $(\;)_{n-1}:n$-$Cat\EE\to (n-1)$-$Cat\EE$. In this way, $(n+1)$-$Cat\EE=Cat^{(\;)_{n-1}}n$-$Cat\EE$. Again, when there is no ambiguity, a $(n+1)$-category will be denoted by the central part $X^{n+1}_{\bullet}$  of the internal category producing it in a fiber of $(\;)_{n-1}:n$-$Cat\EE\to (n-1)$-$Cat\EE$:
$$\xymatrix@=25pt
{
	X^{n+1}_{\bullet}\times_{n-1} X^{n+1}_{\bullet} \ar[rr]|{d^{n+1,1}_{\bullet}} \ar@<2ex>[rr]^{d^{n+1,2}_{\bullet}}\ar@<-2ex>[rr]_{d^{n+1,0}_{\bullet}} && X^{n+1}_{\bullet} \ar@<-2ex>[rr]_{d^{n+1,0}_{\bullet}} \ar@<2ex>[rr]^{d^{n+1,1}_{\bullet}} && X^n_{\bullet} \ar[ll]|{s^{n+1,0}_{\bullet}}  
}
$$
where the left hand side object is a pullback in this fiber. The internal category $	X^n_{\bullet}$ is called the category of $n$-morphisms, while the internal category $X^{n+1}_{\bullet}$ is called the category of $(n+1)$-morphisms or $(n+1)$-cells.  Again for sake of simplicity, we shall denote:\\
1) by $\Pt_{n-1}Cat\EE$ the category  $\Pt_{(\;)_{n-1}}n$-$Cat\EE$ whose objects are the split epimorphisms between internal $n$-categories lying in a fiber of $(\;)_{n-1}$,\\
2) by $\P^C_{n-1}\EE$ the fibration $\P_{(\;)_{n-1}}: \Pt_{n-1}n$-$Cat\EE\to n$-$Cat\EE$ associating with any split epimorphism of this kind its codomain,\\
3) and by $(\GG^C_n,\sigma^C_n,\pi^C_n)$ the monad $(\GG_{(\;)_{n-1}},\sigma_{(\;)_{n-1}},\pi_{(\;)_{n-1}})$ on the category $\Pt_{n-1}n$-$Cat\EE$. Now translating the results of Section \ref{locfib} we get:
\begin{theo}\label{caracnC}
	1) The category $\Alg \GG^C_n$ is the full subcategory $(n+1)_G$-$Cat\EE$ of $(n+1)$-$Cat\EE$ whose objects are the internal groupoids in the fibers of $(\;)_{n-1}: n$-$Cat\EE\to (n-1)$-$Cat\EE$, namely the $(n+1)$-categories with invertible $(n+1)$-cells. The forgetful funtor $(\;)_{n-1}: n$-$Cat\EE\to (n-1)$-$Cat\EE$ associating the $n$-category $X^n_{\bullet}$ with the $(n+1)$-category $X^{n+1}_{\bullet}$  is a fibered reflection.\\ 
    2) The full subcategory of $\GG^C_n$-$Cat(\Pt_{n-1}n$-$Cat\EE)$ whose objects are the $\GG^C_n$-categories  with a $\P^C_{n-1}\EE$-cartesian $0$-leg and an idomorphic $1$-leg is isomorphic to the category of internal categories in the fibers of $(\;)_{n-1}$, namely to $(n+1)$-$Cat\EE$.\\
    3) The inclusion $\GGC^C_n: \Alg\GG^C_n \into (n+1)$-$Cat\EE$ coincides with the following upper horizontal one in the following commutative diagram:
	$$\xymatrix@=20pt
	{
		(n+1)_{G}{\rm -}Cat\EE \ar@{ >->}[r] \ar[d]_{(\;)_n} &	(n+1){\rm -}Cat\EE \ar[d]^-{(\;)_n}  \\
		n-Cat\EE \ar@{=}[r] & n-Cat\EE 
	}
	$$
	where the vertical functors are fibered reflections.
\end{theo}
\proof
Apply Proposition \ref{caracF} and Theorem \ref{mmmainF}.
\endproof

Accordingly, the construction of the following tower of fibered reflections is entirely ruled by the monad $(\GG,\sigma,\pi)$:
$$...\; n{\rm -}Cat\EE \stackrel{(\;)_{n-1}}{\rightarrow} (n-1){\rm -}Cat\EE \; ......\; 2{\rm -}Cat\EE \stackrel{(\;)_1}{\rightarrow} Cat\EE \stackrel{(\;)_0}{\rightarrow} \EE$$
   
\section{The $T_{X_{\bullet}}$-categories}

\subsection{The general case}

In section \ref{T_X} we observed that, when $X_{\bullet}$ is an internal category, the  monad $(T_{X_{\bullet}},\lambda_{X_{\bullet}},\mu_{X_{\bullet}})$ on $\EE/X_0$ is cartesian, and that the algebras of this  monad coincide with the discrete fibrations above $X_{\bullet}$, so that $\AlgT_{X_{\bullet}}=DisF/X_{\bullet}$. We are now going to investigate what are that the $T_{X_{\bullet}}$-categories.
\begin{prop}\label{TTX}
	Given any internal category $X_{\bullet}$ in the category $\EE$, then the category $T_{X_{\bullet}}$-$Cat(\EE/X_0)$ is isomorphic to $Cat\EE/X_{\bullet}$.
\end{prop}
\proof
A pointed $T_{X_{\bullet}}$-graph on an object $g_0:Y_0\to X_0$ of $\EE/X_0$ is given by a diagram of the following kind in $\EE$, where $g_0.d_0^Y=\gamma=d_0.d_1^*(g_0).\bar d_1=d_0.g_1$:
$$\xymatrix@=7pt{
	&&  Y_1 \ar[ddddll]^>>>>>>{\gamma} \ar@<-1ex>[lldd]_{d_0^Y} \ar[ddrr]^{\bar d_1=(d^Y_1,g_1)} &&	   \\
	\\
	Y_0	\ar[dd]_{g_0} \ar[rruu]|{s_0^Y} \ar[rrrr]|{\sigma_0^{g_0}} &&   && d_1^*(Y_0)  \ar@{.>}@<-1,5ex>[llll]_{} \ar[dd]^{d_1^*(g_0)}\\
	&& \\
	X_0  \ar[rrrr]|{s_0} &&  && X_1 \ar@<-1,5ex>@{.>}[llll]_{d_1} \ar@<1,5ex>[llll]^{d_0}
}
$$
satisfying $(d_1^Y,g_1).s_0^Y=\sigma_0^{g_0}=(1_{Y_0},s_0.g_0)$, namely $d_1^Y.s_0^Y=1_{Y_0}$ and $g_1.s_0^Y=s_0.g_0$. Accordingly it is equivalent to a morphism of internal reflexif graphs in $\EE$:
$$\xymatrix@=12pt
{ Y_0 \ar[rr]|{s_0^Y}  \ar[dd]_{g_0}&&
	Y_1    \ar[dd]^{g_1} \ar@<2ex>[ll]^{d_{1}^Y} \ar@<-2ex>[ll]_{d_{0}^Y}\\
	&& \\
	X_0  \ar[rr]|{s_0} &&
	X_1     \ar@<2ex>[ll]^{d_{1}} \ar@<-2ex>[ll]_{d_{0}}
}
$$
We have to build now the pullback of $\bar d_1=(d_1^Y,g_1)$ along $T_{X_{\bullet}}(d_0^Y)$, namely the pullback of $d_1^Y$ along $d_0^Y$ whose domain is denoted $Y_2$:
$$\xymatrix@=7pt{
	Y_1 \ar[dd]_{d_0^Y} \ar@<1,5ex>[rr]|{}  &&	 d_1^*(Y_1) \ar[ll]^{\delta_1^{g_0.d_0}} \ar[dd]^{T_{X_{\bullet}}(g_0)} &&& Y_2 \ar[dd]^{d_0^Y} \ar[lll]_{(d_2^Y,g_1.d_0^Y)}\\
	&&	\\
	Y_0	\ar[dd]_{g_0}  \ar@<1,5ex>[rr]|>>>>>{\sigma_0^{g_0}}  && d_1^*(Y_0)  \ar[ll]^{\delta_1^{g_0}} \ar[dd]^{d_1^*(g_0)} &&& Y_1 \ar[lll]^{(d_1^Y,g_1)}\\
	&& \\
	X_0  \ar[rr]|{s_0} &&  X_1 \ar@<1,5ex>[ll]^{d_1}
}
$$
This induces a map $g_2:Y_2\to X_2$ such that $g_1.d_0^Y=d_0.g_2$ and $g_1.d_2^Y=d_2.g_2$. Accordingly, we get the following diagram where the two central "vertical" triangles commute and where $g_2=d_2^*(g_1).\bar d_2$, with $\bar d_2=(d_2^Y,g_1.d_0^Y)$:

$$\xymatrix@=25pt{
	&&& &&&& Y_2   \ar@{-->}@<1ex>[dll]^{\bar d_2}\ar@{-->}@<-2ex>[dllll]_{d_0^Y} \ar@{.>}[dllll]^{d_1^Y}\\
	&&& Y_1 \ar@{-->}[dl]_>>>>{\bar d_1} \ar[ddl]^<<<<<{g_1} \ar@{-->}[dlll]_{d_0^Y} &&  d_2^*(Y_1) \ar[ddl]^{d_2^*(g_1)} \ar@{-->}[dl]_>>>>{T_{X_{\bullet}}(\bar d_1)} \ar@<1ex>[ll]^>>>>>>>>{\delta_1^{d_0.g_1}}\\
	Y_0  \ar@{ >->}[rr]|{\sigma_0^{g_0}} \ar[d]_{g_0} && d_1^*(Y_0) \ar@<1ex>[ll]^{\delta_1^{g_0}} \ar@<-2pt>[d]_{d_1^*({g_0})} && (d_1.d_2)^*(Y_0) \ar@<-3pt>[d]_{(d_1.d_2)^*(g_0)}\ar@<1ex>[ll]^<<<<<{\delta_2^{g_0}} \ar@<-1ex>[ll]_<<<<<{\delta_1^{g_0}}\\
	X_0   \ar@{ >->}[rr]|{s_0}  && X_1 \ar[d]^{d_0}\ar@<1ex>[ll]^{d_1} && X_2 \ar[d]^{d_0}\ar@<1ex>[ll]^{d_2}\ar@<-1ex>[ll]_{d_1}\\
	&& X_0 && X_1 \ar@<1ex>[ll]^{d_1}\ar@<-1ex>[ll]_{d_0} 
}
$$
The structure of $T_{X_{\bullet}}$-category on $g_0$ is then completed by the data of a map $d_1^Y: Y_2\to Y_1$ in $\EE/X_0$ such that Burroni's Axioms 4, 7, 8 hold. The first part of Axioms 4 is $d_0^Y.d_1^Y=d_0^Y.d_0^Y$ (which implies that $d_1^Y$ is a map in the slice category $\EE/X_0$), while the second part is $\bar d_1.d_1^Y=\mu_{X_{\bullet}}(g_0).T_{X_{\bullet}}(\bar d_1).\bar d_2=\delta_1^{g_0}.T_{X_{\bullet}}(\bar d_1).\bar d_2$. This second part is equivalent to  $d_1^Y.d_1^Y=d_1^Y.d_2^Y$ and $g_1.d_1^Y=d_1.g_2$, which would complete the structure of an internal functor:
$$\xymatrix@=12pt
{ Y_0 \ar[rr]|{s_0^Y}  \ar[dd]_{g_0}&&
	Y_1    \ar[dd]^{g_1} \ar@<2ex>[ll]^{d_{1}^Y} \ar@<-2ex>[ll]_{d_{0}^Y} && Y_2 \ar[dd]^{g_2}  \ar@<2ex>[ll]^{d_{2}^Y}\ar@<-2ex>[ll]_{d_{0}^Y}\ar[ll]|{d_{1}^Y}\\
	&& &&\\
	X_0  \ar[rr]|{s_0} &&
	X_1     \ar@<2ex>[ll]^{d_{1}} \ar@<-2ex>[ll]_{d_{0}} && X_2  \ar@<2ex>[ll]^{d_{2}}\ar@<-2ex>[ll]_{d_{0}}\ar[ll]|{d_{1}}
}
$$
provided that neutrality and associativity of the composition map $d_1^Y: Y_2\to Y_1$ hold, which is straightforward with Axioms 7 and 8.
\endproof
The inclusion $Cat(\EE/X_0)=Cat\EE/\Delta_{X_0}\into T_{X_{\bullet}}$-$Cat(\EE/X_0)=Cat\EE/X_{\bullet}$ is given by composition, in $Cat\EE$, with the inclusion functor $\Delta_{X_0}\into X_{\bullet}$; while its coadjoint $\RR$ (see Proposition \ref{Rad}) is obtained by the pullback in $Cat\EE$ along this inclusion functor. 

\subsection{The $T_{X_{\bullet}}$-groupoids}

\begin{prop}\label{TTXG}
	Given any internal category $X_{\bullet}$ in the category $\EE$, a $T_{X_{\bullet}}$-groupoid is functor above $X_{\bullet}$ whose domain is a groupoid. Accordingly, the category $T_{X_{\bullet}}$-$Grd(\EE/X_0)$ is given by the following pullback:
	$$\xymatrix@=8pt
	{ T_{X_{\bullet}}{\rm -}Grd(\EE/X_0)  \ar@{ >->}[rr]  \ar[dd]_{}&&
		Cat\EE/X_{\bullet}    \ar[dd]^{dom}  \\
		&& &&\\
		Grd\EE  \ar@{ >->}[rr] && Cat\EE 
	}
	$$
\end{prop}
\proof
According to Section \ref{Tgroup}, a $T_{X_{\bullet}}$-category gives rise to a $T_{X_{\bullet}}$-groupoid if and only if the map $d_1^Y:Y_2\to Y_1$ produces the following kernel equivalence relation:
$$\xymatrix@=12pt
{ Y_0 &&
	Y_1     \ar[ll]_{d_{0}^Y} && Y_2   \ar@<-1ex>[ll]_{d_{0}^Y}\ar@<1ex>[ll]^{d_{1}^Y}
}
$$
which is equivalent to the fact that $Y_{\bullet}$ is a groupoid.
\endproof

\subsection{The $T_{X_{\bullet}}$-categories when $X_{\bullet}$ is a groupoid}\label{TTTT_X}

By Section \ref{truc}, we know that:
$$
\xymatrix@C=45pt{
	T_{X_{\bullet}}     &
	T^2_{X_{\bullet}} \ar@{->>}[l]_-{\mu_{X_{\bullet}}} & T^3_{X_{\bullet}}  \ar@<-5pt>[l]_-{\mu_{X_{\bullet}T_{X_{\bullet}}}}
	\ar@<5pt>[l]^-{T_{X_{\bullet}}(\mu_{X_{\bullet}})} }
$$
is a kernel equivalence relation if and only if $X_{\bullet}$ is a groupoid. In this case, by Corollary \ref{TTTTT_X}, any $T_{X_{\bullet}}$-algebra produces a $T_{X_{\bullet}}$-groupoid and we get the following string of inclusions:
$$
\xymatrix@C=25pt {
  \AlgT_{X_{\bullet}} \ar@{ >->}[rr]^<<<<<<<<<<<<{\TXC} \ar@{=}[d] && T_{X_{\bullet}}{\rm -}Grd(\EE/X_0) \ar@{ >->}[r] \ar@{=}[d] & T_{X_{\bullet}}{\rm -}Cat(\EE/X_0) \ar@{=}[d]\\
  DFib/X_{\bullet} \ar@{ >->}[rr] && Grd\EE/X_{\bullet} \ar@{ >->}[r] & Cat\EE/X_{\bullet}
} 
$$

\section{$T$-operads and $T$-multicategories}

About thirty years after Burroni's work \cite{AB} (which was published in french), his ideas have been independantly rediscovered by Leinster \cite{Lein2} and Hermida \cite{H} in the following context.

According to the historical note, p. 63, of Leinster's encyclopedia about operads \cite{Lein}, the notions of operad and multicategory gradually emerged from multiple horizons until they found a name, the first one in May \cite{May} and the second  one in Lambek \cite{Lamb}, before being completely stabilized. Finally, starting with $\EE=Set$ and $(M,\lambda,\mu)$ the free monoid monad which is cartesian as we recalled above, \emph{operads appeared to coincide with $M$-categories with only one object}, while \emph{muticategories} \emph{appeared to coincide with $M$-categories} \cite{Lein2}. Then Leinster introduced the terminolgy $T$-operads and $T$-multicategory for the same notions related to any cartesian monad $(T,\lambda,\mu)$. So, in the cartesian context, $T$-multicategory in the sense of Leinster coincides with $T$-category in the sense of Burroni.

So, given any cartesian monad $(T,\lambda,\mu)$ and following our results, and with respect to our notations related to the inclusion $\bar F^T: \EE \into KlT$, a $T$-multicategory in $\EE$ is nothing but an internal category in $KlT$:
	$$\xymatrix@=18pt{
		X_3 \ar@<-6ex>[rr]_{d_0} \ar@<6ex>@{.>}[rr]^{"\delta_3"} \ar@<-3ex>[rr]|{d_1} \ar@<3ex>[rr]|{d_2} && X_2\ar@<-4ex>[rr]_{d_0} \ar[rr]|{d_1} \ar@<4ex>@{.>}[rr]^{"\delta_2"} \ar[ll]|{s_1}  &&	X_1 \ar@<-2ex>[rr]_{d_0} 
		\ar@<2ex>@{.>}[rr]^{"\delta_1"} \ar@<2ex>[ll]|{s_0} \ar@<-2ex>[ll]|{s_1} && X_0 \ar[ll]|{s_0}
	}
	$$
\textbf{Warning:} Leinster's designation of a $T$-operad in terms of "generalized monoid" could be a bit confusing, because, beyond the undisputable existence of a unit $e$ and of an internal "operation" $m$, a $T$-operad is an actual internal category in $KlT$:
$$\xymatrix@=18pt{
	 X\times_{\delta_1} X \ar@<-4ex>[rr]_{\tau^X_0} \ar[rr]|{m} \ar@<4ex>@{.>}[rr]^{"\delta_2"}   &&	X \ar@<-2ex>[rr]_{\tau_X} 
	\ar@<2ex>@{.>}[rr]^{"\delta_1"} \ar@<2ex>[ll]|{s_0} \ar@<-2ex>[ll]|{s_1} && 1 \ar[ll]|{e}
}
$$
since the object $1$ does not stay a terminal object in $KlT$, unless $T(1)\simeq 1$, and consequently the map $m$ is far from being a classical binary operation.

\subsection{The cartesian monad $(T_{X^T_{\bullet}},\lambda_{X^T_{\bullet}},\mu_{X^T_{\bullet}})$}

Given any cartesian monad $(T,\lambda,\mu)$ on $\EE$ and any $T$-category $X^T_{\bullet}$, Leinster introduced in \cite{Lein2} a notion of algebras associated with them. Indeed,
on the model of Section \ref{T_X}, we get a cartesian monad $(T_{X^T_{\bullet}},\lambda_{X^T_{\bullet}},\mu_{X^T_{\bullet}})$ on the slice category $\EE/X_0$:\\
1) we get a cartesian functor on $\EE/X_0$ since, in the cartesian context, $\EE$ becomes a pullback stable subcategory of $KlT$:
	$$\xymatrix@=20pt{
	Z  \ar@{ >->}[r]|{\sigma_0^h} \ar[d]_{h} & d_1^*(Z) \ar@<1ex>@{.>}[l]^{"\delta_1^h"} \ar@<-2pt>[d]^{d_1^*(h)} & (d_1.d_2)^*(Z) \ar@<-3pt>[d]^{(d_1.d_2)^*(h)}\ar@<1ex>@{.>}[l]^{"\delta_2^h"} \ar@<-1ex>[l]_{\delta_1^h}\\
	X_0   \ar@{ >->}[r]|{s_0}  & X_1 \ar[d]^{d_0}\ar@<1ex>@{.>}[l]^{"\delta_1"} & X_2 \ar[d]^{d_0}\ar@<1ex>@{.>}[l]^{"\delta_2"}\ar@<-1ex>[l]_{d_1}\\
	& X_0 & X_1 \ar@<1ex>@{.>}[l]^{"\delta_1"}\ar[d]^{d_0}\ar@<-1ex>[l]_{d_0} \\
	& & X_0 \\
	h \ar@{ >->}[r]_{\sigma_0^h} & {\;\;T_{X^T_{\bullet}}(h)} & {\;\;\;T^2_{X^T_{\bullet}}(h)} \ar[l]^{\delta_1^h}
}
$$
2) and, in addition, since the diagram $X_0\stackrel{s_0}{\rightarrowtail} X_1 \stackrel{d_1}{\leftarrow} X_2$ lies in $\EE$, so does:\\ $Z\stackrel{\sigma_0^h}{\rightarrowtail} d_1^*(Z) \stackrel{\delta_1^h}{\leftarrow} (d_1.d_2)^*(Z)$; accordingly, this monad is entirely defined inside the slice category $\EE/X_0$. So, in our terms, we get:
\begin{prop}
When $(T,\lambda,\mu)$ is a cartesian monad, the algebras of the monad $(T_{X^T_{\bullet}},\lambda_{X^T_{\bullet}},\mu_{X^T_{\bullet}})$ are the \emph{discrete fibrations} above $X^T_{\bullet}$ in the category $CatKlT$:
$$\xymatrix@=12pt
{ Z_0 \ar[rr]|{s_0}  \ar[dd]_{h_0}&&
	Z_1    \ar[dd]^{h_1} \ar@<2ex>@{.>}[ll]^{"\delta_1"} \ar@<-2ex>[ll]_{d_{0}}\\
	&& \\
	X_0  \ar[rr]|{s_0} &&
	X_1     \ar@<2ex>@{.>}[ll]^{"\delta_1"} \ar@<-2ex>[ll]_{d_{0}}
}
$$
in other words, we get $\AlgT_{X^T_{\bullet}}=DisF(T$-$Cat\EE/X^T_{\bullet})$.
\end{prop}
\proof Straightforward from the classical result on the monad $(T_{X_{\bullet}},\lambda_{X_{\bullet}},\mu_{X_{\bullet}})$\endproof

\subsection{The $T_{X^T_{\bullet}}$-categories}

It remains to characterize the $T_{X^T_{\bullet}}$-categories. On the model of what happens for an ordinary internal category $X_{\bullet}$, we get:

\begin{prop} The $T_{X^T_{\bullet}}$-categories are the $T$-functors above $X^T_{\bullet}$:
$$\xymatrix@=12pt
{ Y_0 \ar[rr]|{s_0}  \ar[dd]_{g_0}&&
	Y_1    \ar[dd]^{g_1} \ar@<2ex>@{.>}[ll]^{"\delta_1"} \ar@<-2ex>[ll]_{d_{0}}\\
	&& \\
	X_0  \ar[rr]|{s_0} &&
	X_1     \ar@<2ex>@{.>}[ll]^{"\delta_1"} \ar@<-2ex>[ll]_{d_{0}}
}
$$
in other words, we get $T_{X^T_{\bullet}}$-$Cat(\EE/X_0)=T$-$Cat\EE/X^T_{\bullet}$.
\end{prop}
\proof Let us follow step by step the proof of Proposition \ref{TTX}. A pointed $T_{X^T_{\bullet}}$-graph on an object $g_0:Y_0\to X_0$ of $\EE/X_0$ is given by a diagram of the following kind in $KlT$, where $g_0.d_0^Y=\gamma=d_0.d_1^*(g_0).\bar d_1=d_0.g_1$ which implies $g_1\in\EE$:
$$\xymatrix@=7pt{
	&&  Y_1 \ar[ddddll]^>>>>>>{\gamma} \ar@<-1ex>[lldd]_{d_0^Y} \ar[ddrr]^{\bar d_1=("\delta^Y_1",g_1)} &&	   \\
	\\
	Y_0	\ar[dd]_{g_0} \ar[rruu]|{s_0^Y} \ar[rrrr]|{\sigma_0^{g_0}} &&   && d_1^*(Y_0)  \ar@{.>}@<-1,5ex>[llll]_{} \ar[dd]^{d_1^*(g_0)}\\
	&& \\
	X_0  \ar[rrrr]|{s_0} &&  && X_1 \ar@<-1,5ex>@{.>}[llll]_{"\delta_1"} \ar@<1,5ex>[llll]^{d_0}
}
$$
satisfying $("\delta^Y_1",g_1).s_0^Y=\sigma_0^{g_0}=(1_{Y_0},s_0.g_0)$, namely $"\delta^Y_1".s_0^Y=1_{Y_0}$ and $g_1.s_0^Y=s_0.g_0$. Accordingly it is equivalent to a morphism of pointed $T$-graphs:
$$\xymatrix@=12pt
{ Y_0 \ar[rr]|{s_0^Y}  \ar[dd]_{g_0}&&
	Y_1    \ar[dd]^{g_1} \ar@<2ex>@{.>}[ll]^{"\delta^Y_1"} \ar@<-2ex>[ll]_{d_{0}^Y}\\
	&& \\
	X_0  \ar[rr]|{s_0} &&
	X_1     \ar@<2ex>@{.>}[ll]^{"\delta_1"} \ar@<-2ex>[ll]_{d_{0}}
}
$$
We have to build now the pullback of $T^T_{X_{\bullet}}(d_0^Y)$ along $\bar d_1=("\delta^Y_1",g_1)$ in $KlT$, which is nothing but the pullback of $d_0^Y$ along $"\delta^Y_1"$ in $KlT$, whose domain is denoted $Y_2$:
$$\xymatrix@=7pt{
	Y_1 \ar[dd]_{d_0^Y} \ar@<1,5ex>[rrr]|{}  &&&	 d_1^*(Y_1) \ar@{.>}[lll]^{"\delta_1^{g_0.d_0}"} \ar[dd]^{T_{X^T_{\bullet}}(g_0)} &&&& Y_2 \ar[dd]^{d_0^Y} \ar[llll]_{("\delta^Y_2",g_1.d_0^Y)}\\
	&&	\\
	Y_0	\ar[dd]_{g_0}  \ar@<1,5ex>[rrr]|>>>>>{\sigma_0^{g_0}} & && d_1^*(Y_0)  \ar@{.>}[lll]^{"\delta_1^{g_0}"} \ar[dd]^{d_1^*(g_0)} &&&& Y_1 \ar[llll]^{("\delta^Y_1",g_1)}\\
	&& \\
	X_0  \ar[rrr]|{s_0} &&&  X_1 \ar@<1,5ex>@{.>}[lll]^{"\delta_1"}
}
$$
This induces a map $g_2:Y_2\to X_2$ in $KlT$ such that $g_1.d_0^Y=d_0.g_2$ (which implies that $g_2$ belongs to $\EE$) and $g_1."\delta^Y_2"="\delta_2".g_2$. Accordingly, we get the following diagram in $KlT$ where the two central "vertical" triangles commute in $\EE$ and where $g_2=d_2^*(g_1).\bar d_2$, with $\bar d_2=("\delta^Y_2",g_1.d_0^Y)$: 

$$\xymatrix@=25pt{
	&&& &&&& Y_2   \ar@<1ex>[dll]^{\bar d_2}\ar@<-2ex>[dllll]_{d_0^Y} \ar@{-->}[dllll]^{d_1^Y}\\
	&&& Y_1 \ar[dl]_>>>>{\bar d_1} \ar[ddl]^<<<<<{g_1} \ar[dlll]_{d_0^Y} &&  d_1^*(Y_1) \ar[ddl]^{d_2^*(g_1)} \ar[dl]_>>>>{T_{X^T_{\bullet}}(\bar d_1)} \ar@<1ex>[ll]^>>>>>>>>{\delta_1^{d_0.g_1}}\\
	Y_0  \ar@{ >->}[rr]|{\sigma_0^{g_0}} \ar[d]_{g_0} && d_1^*(Y_0) \ar@<1ex>@{.>}[ll]^{"\delta_1^{g_0}"} \ar@<-2pt>[d]_{d_1^*({g_0})} && (d_1.d_2)^*(Y_0) \ar@<-3pt>[d]_{(d_1.d_2)^*(g_0)}\ar@<1ex>@{.>}[ll]^<<<<<{"\delta_2^{g_0}"} \ar@<-1ex>[ll]_<<<<<{\delta_1^{g_0}}\\
	X_0   \ar@{ >->}[rr]|{s_0}  && X_1 \ar[d]^{d_0}\ar@<1ex>@{.>}[ll]^{"\delta_1"} && X_2 \ar[d]^{d_0}\ar@<1ex>@{.>}[ll]^{"\delta_2"}\ar@<-1ex>[ll]_{d_1}\\
	&& X_0 && X_1 \ar@<1ex>@{.>}[ll]^{"\delta_1"}\ar@<-1ex>[ll]_{d_0} 
}
$$
The structure of $T_{X^T_{\bullet}}$-category on $g_0$ is then completed by the data of a map $d_1^Y: Y_2\to Y_1$ in $\EE/X_0$ (and thus in $\EE$) such that Burroni's Axioms 4, 7, 8 hold. The first part of Axioms 4 is $d_0^Y.d_1^Y=d_0^Y.d_0^Y$, while the second part is $\bar d_1.d_1^Y=\mu_{X^T_{\bullet}}(g_0).T_{X^T_{\bullet}}(\bar d_1).\bar d_2=\delta_1^{g_0}.T_{X^T_{\bullet}}(\bar d_1).\bar d_2$. This second part is equivalent to  $"\delta^Y_1".d_1^Y="\delta^Y_1"."\delta^Y_2"$ and $g_1.d_1^Y=d_1.g_2$, which would complete the structure of a $T$-functor:
$$\xymatrix@=12pt
{ Y_0 \ar[rr]|{s_0^Y}  \ar[dd]_{g_0}&&
	Y_1    \ar[dd]^{g_1} \ar@<2ex>@{.>}[ll]^{\delta^Y_1"} \ar@<-2ex>[ll]_{d_{0}^Y} && Y_2 \ar[dd]^{g_2}  \ar@<2ex>@{.>}[ll]^{"\delta^Y_2"}\ar@<-2ex>[ll]_{d_{0}^Y}\ar[ll]|{d_{1}^Y}\\
	&& &&\\
	X_0  \ar[rr]|{s_0} &&
	X_1     \ar@<2ex>@{.>}[ll]^{"\delta_1"} \ar@<-2ex>[ll]_{d_{0}} && X_2  \ar@<2ex>@{.>}[ll]^{"\delta_2"}\ar@<-2ex>[ll]_{d_{0}}\ar[ll]|{d_{1}}
}
$$
provided that neutrality and associativity of the composition map $d_1^Y: Y_2\to Y_1$ hold, which is straightforward with Axioms 7 and 8.\endproof
So, the canonical inclusion is the following one: $$\TXTC: \AlgT_{X^T_{\bullet}}=DisF(T{\rm -}Cat\EE/X^T_{\bullet})\into T{\rm -}Cat\EE/X^T_{\bullet}$$

\end{document}